\documentclass[aap,dvips]{arximspdf}

\doi{10.1214/105051606000000835}
\volume{17}
\issue{2}
\pubyear{2007}
\firstpage{537}
\lastpage{571}

\makeatletter
\newcommand{\eqref}[1]{\textup{(\ref{#1})}}
\newtheorem{Lemma}{Lemma}

\newtheorem{Proposition}{Proposition}

\newproclaim{Remark}{Remark}
\newproclaim{Example}{Example}

\def\bX{\mathbf{X}}

\def\bzero{\mathbf{0}}

\def\bx{\mathbf{x}}
\def\bz{\mathbf{z}}

\let\widebar\overline
\makeatother

\begin{document}
\begin{frontmatter}

\title{Limit laws for random vectors with an~extreme~component}
\runtitle{Random Vectors With an Extreme Component}

\begin{aug}
\author[A]{\fnms{Janet E.} \snm{Heffernan}\ead[label=e1]{j.heffernan@lancaster.ac.uk}} and
\author[B]{\fnms{Sidney I.} \snm{Resnick}\corref{}\ead[label=e2]{sir1@cornell.edu}\thanksref{t1}}
\runauthor{J. E. Heffernan and S. I. Resnick}
\affiliation{Lancaster University and Cornell University}
\thankstext{t1}{Supported in part by NSF Grant DMS-03-03493.}
\address[A]{Department of Mathematics and Statistics\\
Lancaster University\\
Lancaster, LA1 4YF\\
United Kingdom\\
\printead{e1}} 
\address[B]{School of Operations Research\\
\quad and Industrial Engineering\\
Cornell University \\
Ithaca, New York 14853\\
USA\\
\printead{e2}}
\end{aug}

\received{\smonth{2} \syear{2005}}
\revised{\smonth{10} \syear{2006}}

%
\begin{abstract}
Models based on assumptions of multivariate regular
variation and hidden regular variation provide ways to describe a
broad range of extremal dependence structures when marginal
distributions are heavy tailed. Multivariate regular variation
provides a rich description of extremal dependence in the case of
asymptotic dependence, but fails to distinguish between exact
independence and asymptotic independence. Hidden regular variation
addresses this problem by requiring components of the random
vector to be simultaneously large but on a smaller scale than the
scale for the marginal distributions. In doing so, hidden regular
variation typically restricts attention to that part of the
probability space where all variables are simultaneously large.
However, since under asymptotic independence the largest values do
not occur in the same observation, the region where variables are
simultaneously large may not be of primary interest. A different
philosophy was offered in the paper of Heffernan and Tawn
[\textit{J.~R.~Stat. Soc. Ser.~B Stat. Methodol.} \textbf{66} (2004)
\mbox{497--546}]
which allows examination of distributional tails other than the
joint tail. This approach used an asymptotic argument which
conditions on one component of the random vector and finds the
limiting conditional distribution of the remaining components as
the conditioning variable becomes large. In this paper, we
provide a thorough mathematical examination of the limiting
arguments building on the orientation of
\mbox{Heffernan} and Tawn [\textit{J.~R.~Stat. Soc. Ser.~B Stat.
Methodol.} \textbf{66} (2004) 497--546]. We examine the conditions required
for the assumptions made by the conditioning approach to hold, and
highlight simililarities and differences between the new and
established methods.
\end{abstract}

%

\end{frontmatter}

\section{Introduction}\label{intro}

Extreme value theory motivates statistical models for the tails of
{multivariate} probability distributions. All such theory relies
on some form of asymptotic argument; it is this limiting argument
which forces us into the distributional tails and allows the
examination of the extremal behavior of random vectors.

The first such arguments relied upon limiting behavior imposed by
considering componentwise maxima of random vectors
\cite{dehaanresnick1977,dehaan1985,pickands1981,resnick1987}. This
approach was extended by Coles and Tawn
\cite{colestawn1991,colestawn1994}, de Haan and de Ronde
\cite{dehaanderonde1998} in a multivariate analogue of the
one-dimensional threshold methods of Davison and
Smith~\cite{smith1989}, Smith~\cite{davisonsmith1990}. The methods provide a
rich class of models to describe asymptotic dependence but cannot
distinguish between asymptotic independence and exact independence. In
response to this weakness, theory and models offering a richer
description of asymptotic independence behavior have been developed by
\mbox{Heffernan} and Resnick \cite{heffernanresnick2004}, Ledford and Tawn
\cite{ledfordtawn1996,ledfordtawn1997,ledfordtawn1998}, Maulik and
Resnick \cite{maulikresnick2003b} and Resnick~\cite{resnick2002a}. The
assumptions underlying this broader class of models have been termed
\textit{hidden regular variation} which elaborates the concept of the
\textit{coefficient of tail dependence}.

Models based on assumptions of multivariate regular variation and
hidden regular variation have a common reliance on limiting procedures
in which all vector components are scaled by functions increasing to
infinity. In the case of asymptotic dependence, reliance only on
multivariate regular variation is sufficient since in this case the
largest values of the components of the random vector tend to occur
together. However, models based on multivariate regular variation fail
to distinguish between asymptotic independence and exact independence
and as such provide an inadequate description of dependence within the
asymptotic independence class. Hidden regular variation attempts to
repair this defect by allowing a different scale function which gives
nontrivial limit behavior when vector components are simultaneously
large. Although the hidden regular variation as typically formulated
provides a more satisfactory description of the joint tail of the
distribution for asymptotically independent variables, this approach
still has practical limitations in applications where interest is in
tail regions other than the joint tail. These other tail regions are of
practical significance since under asymptotic independence, the largest
values of the components of the random vector tend not to occur in the
same observation.

The philosophy of examining distributional tails in which one or more
but not necessarily all of the vector components are simultaneously
large was explained in \cite{heffernantawn2004}. They focused on a
single variable being large by conditioning on one component of the
random vector and finding the limiting conditional distribution of the
remaining components as the conditioning variable becomes large.
Simulation studies in~\cite{heffernantawn2004} suggested that this
alternative approach is useful in accurately describing a range of
qualitatively different dependence structures including asymptotic
dependence, asymptotic independence and negative dependence. The
approach is flexible and readily applicable to general $d$-dimensional
distributions. However, this new basis for modeling multivariate
extremes was criticized in the discussion to the paper as lacking a
rigorous theoretical underpinning. The discussion highlighted the need
for further work to clarify how the approach extends and/or differs
from established methodologies which rely on multivariate regular
variation and hidden regular variation.

In this paper, we use the philosophy of \mbox{Heffernan} and Tawn
\cite{heffernantawn2004} and offer a mathematical framework for a theory
of conditional distributions given a component is large. We have
changed the formulation of \mbox{Heffernan} and Tawn
\cite{heffernantawn2004} for two reasons. First, it is difficult to
construct an asymptotic theory based on regular conditional
distributions which are readily manageable only for the case in which
smooth densities are assumed and secondly our formulation readily
allows for connections to classical multivariate extreme value theory
and regular variation.

\subsection{Content of the paper}\label{subsec:content}
Here are more details about the content of the paper.
We consider the distribution of a bivariate random vector $(X,Y)$ on
$\mathbb{R}^2$ under the condition that $Y$ is large. Generalizations
could be made to the case of a~\mbox{$(d+1)$-}dimensional vector
\[
(\bX,Y):=\bigl(X^{(1)},\dots,X^{(d)}, Y\bigr)
\]
where we seek conditional
limits of $\bX$ given $Y$ is large. However, we leave such
generalizations to subsequent investigations. We assume the
distribution function $F$ of $Y$ is in a \textit{domain of attraction}
of an extreme value distribution $G_\gamma(x) $, written $F \in
D(G_\gamma) $. This means there exist functions $a(t)>0, b(t) \in
\mathbb{R}$, such that,
%
\begin{equation}\label{eqn:doa}
F^t\bigl(a(t)y +b(t)\bigr) \to G_\gamma(y)\qquad (t\to\infty),
\end{equation}
weakly, where
%
\begin{equation}\label{eqn:evdistn}
G_\gamma(y)=\exp\{ -(1+\gamma y)^{-1/\gamma} \},\qquad  1+\gamma y >0,
\gamma\in\mathbb{R},
\end{equation}
and the expression on the right is interpreted as {$e^{-e^{-y}}$} if
$\gamma=0.$ See, for example,
\cite{coles2001,embrechtskluppelbergmikosch1997,dehaan1970,reissthomas2001,resnick1987}.
We can and do assume
\[
b(t)= \biggl( \frac{1}{1-F(\cdot)} \biggr)^\leftarrow(t),
\]
where for a nondecreasing function $U$ we define the left continuous
inverse
\[
U^\leftarrow(t)=\inf\{y\dvtx U(y)\geq t\}.
\]
Setting $\widebar F =1-F$, we have relation \eqref{eqn:doa} is
equivalent to
%
\begin{equation}\label{eqn:doa'}
t\widebar F\bigl(a(t)y + b(t)\bigr) \to(1+\gamma y)^{-1/\gamma},\qquad
1+\gamma y
>0,
\end{equation}
or taking inverses
%
\begin{equation}\label{eqn:doa2prime}
\frac{b(tx)-b(t)}{a(t)} \to\frac{x^\gamma-1}{\gamma},\qquad  x>0.
\end{equation}

For convenience we write $\mathbb{E}_\gamma:=\{y\in\mathbb{R}\dvtx 1+\gamma
y>0\}$. When considering vague convergence, it is convenient to close
the interval $\{y\in\mathbb{R}\dvtx 1+\gamma y>0\}$ on the right and
denote by $\widebar \mathbb{E}_\gamma$ this closure. So, for instance,
$\widebar \mathbb{E}_0=(-\infty,\infty]$.

In Section~\ref{sec:basicDefResults}, we
explore the implications of assuming the existence of:
\begin{enumerate}
\item\label{item:scaling}
Scaling function $a(\cdot)>0,$ and centering function
$b(\cdot) \in\mathbb{R}$ so that
\eqref{eqn:doa} holds for $F(x)=P[Y\leq x]$;

\item\label{item:centering} Scaling function $\alpha(\cdot)>0$, and
centering function $\beta(\cdot)\in\mathbb{R}$ and a nonnull Radon
measure $\mu$ on Borel subsets of
$[-\infty,\infty]\times(-\infty,\infty]$, such that for each fixed~$y
\in\mathbb{E}_\gamma$,
\begin{longlist}[(a)]
\item[(a)] $\mu ( [-\infty,x]\times(y,\infty] )$ is not a degenerate
distribution function in $x$,

\item[(b)] $\mu([-\infty,x]\times(y,\infty])<\infty$,

\item[(c)] and
%
\begin{equation}\label{eqn:basicConv}
tP \biggl[\frac{X-\beta(t)}{\alpha(t)} \leq x, \frac{Y-b(t)}{a (t)} >y
\biggr] \to\mu \bigl( [-\infty,x]\times(y,\infty] \bigr),
\end{equation}
at continuity points $(x,y)$ of the limit.
\end{longlist}
\end{enumerate}
If we interpret \eqref{eqn:basicConv} as vague convergence (cf.
Section~\ref{subsubsec:vagueConv}) in $M_+([-\infty,\infty]\times
\widebar \mathbb{E}_\gamma),$ the Radon measures on $[-\infty,\infty]\times
\widebar \mathbb{E}_\gamma$, then in fact \eqref{eqn:basicConv} implies $F \in
D(G_\gamma)$ for some $\gamma\in\mathbb{R}$. Also, we will see that
\eqref{eqn:basicConv} is equivalent to assuming the existence of the
conditional limiting distribution of the scaled and centered $X$
variable given $Y$ is extreme:
%
\begin{equation}\label{eqn:conditVersion}
P \biggl[ \frac{X-\beta\circ b^\leftarrow(t) }{\alpha\circ
b^\leftarrow(t)} \leq x \Big| Y>t \biggr] \to\mu\bigl(
[-\infty,x]\times(0,\infty] \bigr),
\end{equation}
as $t$ converges to the right end point of $F$. This observation
motivates our focusing on the convergence \eqref{eqn:basicConv}.

Thus we make a different assumption from that of \mbox{Heffernan} and Tawn
\cite{heffernantawn2004}, in that in~(\ref{eqn:conditVersion}) we
condition on the event $Y>t$ rather than $Y=t$ as
in~\cite{heffernantawn2004} which requires regular conditional
distributions which are only defined up to almost everywhere
equivalence. Our formulation also has a natural connection with extreme
value theory as it implies $Y$ is in a domain of attraction. In cases
where densities exist, the two formulations are similar. See
Section~\ref{subsec:densities}.

Having established conditions for the existence of a limit
in~\eqref{eqn:basicConv}, in Section~\ref{sec:charLimits} we
characterize the class of attainable limiting measures. These measures
are found to be either product measures or to have a spectral form
after a standardization procedure and then transformation to polar
coordinates. The standardization renders \eqref{eqn:basicConv} into a
standard multivariate regular variation condition on the cone
$[0,\infty]\times(0,\infty]$ and puts us in familiar territory.
Relating \eqref{eqn:basicConv} to standard multivariate regular
variation allows us to identify the class of possible limit measures
\cite{resnickbook2006,resnick1987,resnick2004b}.

Section~\ref{sec:randomNorming} is motivated by the \mbox{Heffernan} and Tawn
\cite{heffernantawn2004} approach. Instead of normalizing $X$ by
deterministic functions of the threshold $t$, we normalize by functions
of the precise value of $Y$ occurring with $X$. This leads to a product
limit form in all cases.

In Section~\ref{sect:ConnectionToAsympIndpce}, we highlight connections
between assumption~\eqref{eqn:basicConv} and standard assumptions of
multivariate regular variation and hidden regular variation, and in
particular show that under multivariate regular
variation,~\eqref{eqn:basicConv} assumes something additional beyond
multivariate regular variation only in the presence of asymptotic
independence.

Section~\ref{sec:egs} illustrates our results with a range of examples.
Of particular interest is the bivariate Normal example which shows a
transformation of $X$ for which the limit~\eqref{eqn:basicConv} does
not exist. This leads to Section~\ref{sec:changeCoord}, in which we
explore how flexible one can be in the choice of measurement units in
which to record $X$ such that the limit measure
in~\eqref{eqn:basicConv} does exist. Our results suggest how to
construct change of variable functions which will give such a limit.

Section~\ref{sect:discuss} returns in more detail to the modeling
assumptions made by \mbox{Heffernan} and Tawn \cite{heffernantawn2004} which
motivated the work of this paper, and discusses the implications of the
new results for their conditional approach to modeling multivariate
extreme values.

\subsection{Symbol and concept glossary}\label{subsec:glossary}
The \hyperref[sec:appendices]{Appendix} contains several appendices
reviewing and referencing needed background. We merely list here some
concepts and symbols; explanations and references in the appendices can
be consulted as needed.
\[
\begin{tabular}{@{}l@{\quad}p{110mm}@{}}
vectors & Bold lower case is reserved for deterministic vectors and
bold upper case is reserved for random vectors. Relations are interpreted
componentwise. See Section~\ref{subsubsec:vecNotation}.
\\
$\mathbb{E}$ & A nice subset of compactified finite dimensional Euclidean space.
\\
$M_+(\mathbb{E})$ & The class of Radon measures on Borel subsets of $\mathbb{E}$.
\\
$U^\leftarrow$ & The left continuous inverse of the nondecreasing function $U$.
\\
$RV_\rho$ & The class of regularly varying functions with index $\rho$ defined
in \eqref{eqn:defRV}.
\\
$\Pi$ & The function class $\Pi$ reviewed in Section~\ref{subsubsec:fctClass} along with
subclasses $\Pi_+(a (\cdot) )$ and $\Pi_- (a(\cdot) )$ and auxiliary function~$a(\cdot)$.
\\
$\Gamma$ & The function class $\Gamma$ reviewed in Section~\ref{subsubsec:fctClass} along with $\Gamma(f)$
and auxiliary function~$f$.
\\
$\stackrel{v}{\to}$ & Vague convergence of measures; see Section~\ref{subsubsec:vagueConv}.
\\
$G_\gamma$ & An extreme value distribution given by \eqref{eqn:evdistn} in the Von Mises
parameterization.\\
\end{tabular}
\]

\[
\begin{tabular}{@{}l@{\quad}p{110mm}@{}}
$\mathbb{E}_\gamma$& $\{x\dvtx 1+\gamma x>0\}$.
\\
$\widebar \mathbb{E}_\gamma$ & The closure on the right of the interval $\mathbb{E}_\gamma$.
\\
$D(G_\gamma)$ & The domain of attraction of the extreme value distribution $G_\gamma$.
This is the set of $F$'s satisfying \eqref{eqn:doa}.
Note for $\gamma>0$, $F \in D(G_\gamma)$
is equivalent to $1-F \in RV_{1/\gamma}$.
\end{tabular}
\]

\section{Basic results}\label{sec:basicDefResults}
In this section we give some implications of
\eqref{eqn:basicConv} and the assumptions (\ref{item:scaling}),
(\ref{item:centering}) given in
Section~\ref{intro}.

\subsection{Standardization of $Y$}
Without loss of generality, we may
assume $Y$ is heavy tailed and $F\in D(G_1).$ The usual
standardization procedure in
extreme value theory (e.g., \cite{resnick1987}, Chapter~5,
\cite{dehaanferreira2006}, Chapter~6.1.2, \cite{resnickbook2006},
Section~6.5.6) means
that \eqref{eqn:doa} implies for $x>0$, as $t \to\infty$,
\begin{eqnarray*}
tP\biggl[\frac{b^\leftarrow(Y)}{t}>x\biggr]&=&
tP\biggl[\frac{Y-b(t)}{a(t)}>\frac{b(tx)-b(t)}{a(t)}\biggr]
\\
&\hspace*{3pt} \to& \biggl( 1+\gamma\frac{(x^\gamma-1)}{\gamma} \biggr)^{-1/\gamma} =x^{-1}.
\end{eqnarray*}
Note if the distribution $F$ of $Y$
is continuous, $b^\leftarrow(Y)$ has a Pareto
distribution and, in any case, $b^\leftarrow(Y)$ will always have a
distribution tail which is asymptotically Pareto. For $y>0$,
\eqref{eqn:basicConv} and \eqref{eqn:doa2prime} imply
%
\begin{eqnarray}\label{eqn:standard form}
\nonumber && tP \biggl[
\frac{X-\beta(t)}{\alpha(t)} \leq x
,\frac{b^\leftarrow(Y)}{t} >y \biggr]
\\
&&\qquad = tP \biggl[\frac{X-\beta(t)}{\alpha(t)} \leq x,
\frac{Y-b(t)}{a (t)} >\frac{b(ty)-b(t)}{a(t)}
\biggr]
\\
\nonumber &&\qquad \to
\cases{
\mu \biggl([-\infty,x]\times\biggl(\dfrac{y^\gamma-1}{\gamma},\infty\biggr]
\biggr),&\quad if $\gamma\neq0$,
\cr
\mu \bigl([-\infty,x]\times(\log y,\infty]
\bigr),&\quad  if $\gamma= 0$.}
\end{eqnarray}
So at the expense of replacing $Y$ by $b^\leftarrow(Y)$,
theoretical development proceeds
without loss of generality by replacing the conditions around
\eqref{eqn:basicConv} with
%
\begin{eqnarray}\label{eqn:basicConvStandardForm}
\hspace*{4mm} \cases{
\mu \bigl([-\infty,x]\times(y,\infty] \bigr) \mbox{ is not a degenerate distribution function in $x$,}
\cr\noalign{}
\qquad \mbox{for each $y>0$},
\cr\noalign{}
P[Y \leq t] \in D(G_1),\qquad \displaystyle{\lim_{t \to\infty}tP[Y>t]=1},
\cr\noalign{}
tP \biggl[\dfrac{X-\beta(t)}{\alpha(t)} \leq x, \dfrac{Y}{t} >y \biggr]
\to\mu\bigl( [-\infty,x]\times(y,\infty] \bigr),
\cr\noalign{}
\qquad x\in \mathbb{R}, y>0, \mbox{ at continuity points $(x,y)$ of the
limit.}}
\end{eqnarray}
We refer to \eqref{eqn:basicConvStandardForm} as the basic
convergence with the $Y$-variable \textit{standardized}.

\begin{Remark}\label{rem:changeCoords}
The argument leading to \eqref{eqn:basicConvStandardForm} shows that we
are free to change the marginal distribution of the $Y$-variable
without disturbing the conditional convergence
\eqref{eqn:conditVersion}. We will see in Section~\ref{sec:egs}, that
this is not always possible for the $X$-variable.
\end{Remark}

We reiterate the connection with conditional modeling when
\eqref{eqn:basicConvStandardForm} is assumed. For $x$
which are continuity points of $H(x):=\mu ( [-\infty,x]\times
(1,\infty] )$,
%
\begin{eqnarray}\label{eqn:conditModelfromConv}
\nonumber H_t \bigl(\alpha(t)x+\beta(t) \bigr)
&:=& P\biggl[\frac{X-\beta(t)}{\alpha(t)} \leq x \Big|Y>t\biggr]
\\
&\hspace*{3pt} =& \frac{P[(X-\beta(t))/\alpha(t) \leq x ,Y>t]}
{P[Y>t]}
\\
\nonumber  &\hspace*{3pt} \sim& tP\biggl[\frac{X-\beta(t)}{\alpha(t)} \leq x ,\frac{Y}{t}>1\biggr]
\\
\nonumber &\hspace*{3pt} \to&
\mu\bigl([-\infty,x]\times(1,\infty] \bigr)=:H(x).
\end{eqnarray}
Interpreting \eqref{eqn:basicConvStandardForm} as vague convergence
on $M_+([-\infty,\infty]\times(0,\infty])$, we obtain from marginal
convergence that
\[
H(\infty)=\mu\bigl([-\infty,\infty]\times(1,\infty]\bigr)=1.
\]

\subsection{Properties of the functions $\alpha(\cdot) $ and $\beta(\cdot)$}
The following is an initial attempt to understand the properties of
the functions $\alpha(\cdot) $ and $\beta(\cdot).$

\begin{Proposition}\label{prop:howVary}
Suppose $(X,Y)$ satisfy the standard form condition
\eqref{eqn:basicConvStandardForm}. Then there exist two functions
$\psi_1(\cdot)$, $\psi_2(\cdot),$ such that for all $c>0$,
%
\begin{equation}\label{eqn:varyAlpha}
\lim_{t\to\infty}\frac{\alpha(tc)}{\alpha(t)} = \psi_1(c)
\end{equation}
and
%
\begin{equation}\label{eqn:varyBeta}
\lim_{t\to\infty}\frac{\beta(tc)-\beta(t)}{\alpha(t)} = \psi_2(c).
\end{equation}

The convergence in \eqref{eqn:varyAlpha} and \eqref{eqn:varyBeta} is
uniform on compact subsets of $(0,\infty)$.
\end{Proposition}

\begin{pf}
Pick $c>0$. For all but an at most countable set $\Lambda$ of
$x$-values, $(x,1)$ and $(x,c^{-1})$ are continuity points of $\mu$.
For $x \in\Lambda^c$, on the one hand we have
\eqref{eqn:conditModelfromConv} and on the
other we have
%
\begin{eqnarray}\label{eqn:defHc}
\nonumber && \lim_{t \to\infty}
P \biggl[
\frac{X-\beta(tc)}{\alpha(tc)} \leq x\bigg|\frac{Y}{t} >1
\biggr]
\\
&&\qquad =\lim_{t \to\infty}
tP \biggl[\frac{X-\beta(tc)}{\alpha(tc)} \leq x, \frac{Y}{t}
>1 \biggr]
\nonumber
\\[-8pt]
\\[-8pt]
\nonumber
&&\qquad = \lim_{t \to\infty}
\frac{tc}{c}P \biggl[
\frac{X-\beta(tc)}{\alpha(tc)} \leq x,\frac{Y}{tc} >c^{-1}
\biggr]
\\
\nonumber &&\qquad = \frac{\mu ( [-\infty,x]\times(c^{-1},\infty] )
}{c}=:H^{(c)}(x).
\end{eqnarray}

Thus the family $\{H_t\}$
converges with two different normalizations:
\[
H_t\bigl(\alpha(t)x+\beta(t)\bigr) \to H(x),\qquad
H_t\bigl(\alpha(tc)x+\beta(tc)\bigr) \to H^{(c)}(x).
\]
The convergence to types theorem (see, e.g.,
\cite{feller1971} or \cite{resnick1998book}, page~275) implies that
\eqref{eqn:varyAlpha} and \eqref{eqn:varyBeta} hold and also
%
\begin{equation}\label{eqn:fromConvToTypes}
H^{(c)}(x) =H\bigl(\psi_1(c)x + \psi_2(c)\bigr).
\end{equation}

To prove local uniform convergence in \eqref{eqn:varyAlpha} and
\eqref{eqn:varyBeta}, replace $c>0 $ in the argument with $c(t)$ where
$c(t) \to c\in(0,\infty)$. Then \eqref{eqn:varyAlpha} and
\eqref{eqn:varyBeta} still hold and since $\psi_1, \psi_2$ are
continuous (see next paragraph), the result follows from
\textit{continuous convergence}. See \cite{resnick1987}, page~2, or~\cite{kuratowski1966}.
\end{pf}

From
\eqref{eqn:varyAlpha}, we have that $\alpha(\cdot) $ is regularly
varying with some index $\rho\in\mathbb{R}$, written $\alpha\in
RV_\rho$, so that $\psi_1(x)=x^{\rho}$. (See \cite{resnick1987}, page~14,
\cite{binghamgoldieteugels1987,feller1971,gelukdehaan1987,dehaan1970,seneta1976}.)
The function $\psi_2(x)$ may be identically zero.
However, if it is not, then from \cite{gelukdehaan1987}, page~16,
we have
%
\begin{equation}\label{eqn:psi2}
\psi_2(x)=
\cases{
k (x^\rho-1)/\rho,&\quad  if $\rho\neq0, x>0$,
\cr\noalign{}
k\log x,&\quad  if $\rho=0, x>0$,}
\end{equation}
for $k \neq0$. Also, there is more detailed information:
\begin{longlist}[(iii)]
\item[(i)] If $\rho>0$, then $\beta(\cdot) \in RV_\rho$ and $\beta(t)
\sim\frac1\rho\alpha(t).$
So it is enough to scale $X$ in~\eqref{eqn:basicConvStandardForm} with a consequent location change in
the $x$-variable for $\mu$.

\item[(ii)] If $\rho=0$, then $\beta(\cdot) \in\Pi(\alpha)$ and
$\alpha\in RV_0$. So $\alpha$ is the auxiliary function of the
$\Pi$-function $\beta$.

\item[(iii)] If $\rho<0$, then $\beta(\infty)= \lim_{t \to\infty}
\beta(t) $ exists finite and
\[
\beta(\infty) -\beta(t) \in
RV_\rho;\qquad  \bigl(\beta(\infty)-\beta(t) \bigr)
\sim\frac{1}{|\rho|}\alpha (t).
\]
\end{longlist}

Case (iii) can be reduced to case (i) by a change of variable.
From case (iii) of~\eqref{eqn:basicConvStandardForm} we get
\[
tP \biggl[ \frac{X-\beta(\infty) +[\beta
(\infty) -\beta(t)]} {
|\rho|
(\beta(\infty)-\beta(t))
} \leq x,
\frac{Y}{t}>y \biggr] \to\mu\bigl( [-\infty,x]\times(y,\infty]
\bigr).
\]
Write
%
\begin{equation}\label{eqn:negRhoToPosRho}
\widetilde X:= \frac{1} {X-\beta(\infty)}, \qquad\tilde\beta (t):=\frac{1} {
{|\rho|} (\beta(\infty) -\beta(t) )},
\end{equation}
so that
%
\begin{eqnarray}\label{eqn:transformedCase}
\nonumber && tP \biggl[ \frac{\widetilde{ X}}{\tilde\beta(t)} \leq x, \frac{Y}{t}>y \biggr]
\\
&&\qquad =
tP \biggl[ \frac{ X- \beta(\infty)} { |\rho| (\beta(\infty)-\beta(t))
}\geq\frac1x,
\frac{Y}{t}>y \biggr]
\nonumber
\\[-8pt]
\\[-8pt]
\nonumber
&&\qquad = tP \biggl[
\frac{ X- \beta(\infty)}
{ |\rho| (\beta(\infty)-\beta(t))
}+\frac{ 1}{{|\rho|} }
\geq\frac1x +\frac{ 1}{{|\rho|}},
\frac{Y}{t}>y \biggr]
\\
\nonumber &&\qquad \to \mu\biggl(
\biggl[\frac1x +\frac{ 1}{{|\rho|}},\infty\biggr]\times(y,\infty]
\biggr)=:\tilde\mu \bigl( [-\infty,x]\times
(y,\infty] \bigr).
\end{eqnarray}
Since case (iii) can be reduced to case (i), it does not
need separate theoretical attention.

\subsection{Conditions for the limit $\mu$ to be a product
measure}\label{subsubsec:whenMuProd}
It turns out that $\mu$ being a product measure is equivalent to
$\psi_1\equiv1$ and $\psi_2\equiv0$.

\begin{Proposition}\label{prop:whenMuProd}
We have $\mu=H\times\nu_1$, where $\nu_1((y,\infty])=y^{-1},
y>0$ (i.e., $\mu( [-\infty,x]\times(y,\infty] )=H(x)y^{-1}$),
iff for all $c>0$,
%
\begin{equation}\label{eqn:whenProd}
{\psi_1(c)=}\lim_{t \to\infty}\frac{\alpha(tc)}{\alpha(t)}=1,\qquad
{\psi_2(c)=}\lim_{t \to\infty} \frac{\beta(tc)-\beta(t)}{\alpha(t)}=0.
\end{equation}
\end{Proposition}

\begin{pf}
Given that $\mu$ is a product, we have from
\eqref{eqn:conditModelfromConv} and \eqref{eqn:defHc}, that
$H^{(c)}(x)=H(x)$. Hence \eqref{eqn:whenProd} follows from the
convergence to types theorem. Conversely, if \eqref{eqn:whenProd}
holds, $H^{(c)}(x)=H(x)$ and from \eqref{eqn:defHc} we have, for all
$c>0$, $\mu( [-\infty, x]\times(c^{-1},\infty] )=cH(x).$ So for all
$y>0$, $\mu( [-\infty,x]\times(y,\infty] )=H(x)y^{-1}.$
\end{pf}

\begin{Remark}\label{rem:psi2ZeroButNotpsi1}
What if $\psi_2\equiv0 $ but $\psi_1 \not\equiv1$? Then $\alpha \in
RV_\rho$ for some $\rho\in\mathbb{R}$, $\rho\not= 0$ and
$\psi_1(c)=c^\rho$, for $c>0$. The reasoning in the previous proof
shows that $\mu$ has the form
%
\begin{equation}\label{eqn:sortOfProduct}
\mu\bigl([-\infty,x]\times(y,\infty]\bigr)=y^{-1}H(x/y^\rho),
\end{equation}
for $x\in\mathbb{R}$, and $y>0$ and where $H$ is a proper nondegenerate
probability distribution.
\end{Remark}

\subsection{When the $X$-variable can be
standardized}\label{subsec:canWeStandardizeX} Standardization is the
process of transforming variables so that their distributions have
regularly varying tails \textit{in standard form}. See
\cite{resnick1987}, Chapter~5, \cite{dehaanferreira2006},
Chapter~6.1.2, \cite{resnickbook2006}, Section 6.5.6. Once standard
form regular variation is achieved, limit measures have a scaling
property and characterization of these limits becomes possible. We know
we can standardize the $Y$ variable. What about the $X$ variable?

It is possible to standardize the $X$-variable if $\beta(t) \geq0$ and
$\psi_2(\cdot)$ in \eqref{eqn:varyBeta} is not constant and
$\beta^\leftarrow$ is nondecreasing on the range of $X$ since in this
case we have for $x>0$,
%
\begin{eqnarray}\label{eqn:standardizeTheSucker}
tP \biggl[ \frac{\beta^\leftarrow(X)}{t} \leq x,\frac Yt >y \biggr] &=&
tP \biggl[ \frac{X-\beta(t)}{\alpha(t)} \leq\frac{\beta(tx) -\beta
(t)}{\alpha(t)},\frac Yt >y \biggr] \nonumber
\\[-8pt]
\\[-8pt]
\nonumber &\hspace*{3pt} \to& \mu \bigl( [-\infty, \psi_2(x)]\times(y,\infty] \bigr),
\end{eqnarray}
at continuity points of the limit. We emphasize there are important
cases where $\psi_2(x)$ is identically zero and thefore where $X$
cannot be standardized by the procedure in
\eqref{eqn:standardizeTheSucker}; see Section~\ref{subsec:bivarNormal}.

Standardization is also possible if $\psi_2 \equiv0$, provided $X>0$
and $\psi_1 \not\equiv1$; that is if $\alpha(\cdot) \in RV_\rho$ with
$\rho\not=0$. If $\rho>0$, then \cite{binghamgoldieteugels1987},
Theorem 3.1.12a,~c, page 136, gives \mbox{$\beta(t)/\alpha(t) \to0$} and by the
convergence to types theorem \eqref{eqn:basicConvStandardForm} can be
rewritten as
\[
tP\biggl[\frac{X}{\alpha(t)} \leq x, \frac Yt >y\biggr] \to
\mu\bigl([0,x]\times(y, \infty]\bigr),\qquad  x>0,y>0.
\]
Therefore, supposing without loss of generality that $\alpha(\cdot)$ is
strictly increasing and continuous (e.g., \cite{seneta1976}), we have
\begin{eqnarray*}
tP\biggl[\frac{\alpha^\leftarrow(X)}{t} \leq x, \frac Yt >y\biggr] &=&
tP\biggl[\frac{ X}{\alpha(t)} \leq\frac{\alpha(t x)}{\alpha(t)}, \frac
Yt >y\biggr]
\\
&\hspace*{3pt} \to& \mu\bigl( (0,x^\rho]\times(y,\infty] \bigr)
\end{eqnarray*}
and $(\alpha^\leftarrow(X),Y) $ are the standardized variables. If
$\rho<0$, \cite{binghamgoldieteugels1987}, Theorem~3.1.10a,~c,
page~134, implies $\beta(\infty):=\lim_{t\to\infty}\beta(t) $ exists finite and $
(\beta(\infty) -\beta(t) ) /\break \alpha(t) \to0.$ Therefore, if we suppose
$P[X~\leq~\beta(\infty)]=1,$ we have for $x>0$,
\begin{eqnarray*}
&& \lim_{t\to\infty} tP \biggl[ \frac{1/(\beta(\infty)-X)}{1/\alpha(t)}
\leq x,\frac Yt>y \biggr]
\\
&&\qquad   =\lim_{t\to\infty}t
P\biggl[\frac{\beta(\infty)-X}{\alpha(t)}\geq x^{-1} ,
\frac Yt>y\biggr]
\\
&&\qquad =\lim_{t\to\infty} tP\biggl[\frac{\beta(\infty) -X
-(\beta(\infty)-\beta(t))}{\alpha(t)}\geq x^{-1} , \frac Yt>y\biggr]
\\
&&\qquad =\lim_{t\to\infty} tP\biggl[\frac{X-\beta(t)}{\alpha(t)}
\leq-x^{-1} , \frac Yt>y\biggr]
\\
&&\qquad =\mu\bigl([-\infty,-x^{-1}]\times(y,\infty]\bigr),
\end{eqnarray*}
and the variables $((\beta(\infty) -X)^{-1},Y)$ can be standardized
according to the recipe for the $\rho>0$ case.

\subsubsection{When $\beta(t)$ is monotone}\label{subsubsec:monotone}
The standardization of the $X$ variable in
\eqref{eqn:standardizeTheSucker} begs the question of when $\beta $ is
monotone. Consider the case where $\psi_2 \not\equiv0$ and $\psi_2$ is
given by \eqref{eqn:psi2} and indexed by $\rho\in \mathbb{R}$. For
discussing when $\beta(t)$ is monotone, it is important to remember
that $\beta(\cdot)$ is only determined up to the asymptotic equivalence
given by the convergence to types theorem.

Consider the following cases.
\begin{longlist}[3.]
\item[1.] $\rho>0$: For this case, we have $\beta\in RV_\rho$ and there
exists $\tilde\beta(t) \in RV_\rho$ such that $\tilde\beta (\cdot)$ is
continuous, strictly increasing to $\infty$ with $\beta
\sim\tilde\beta$. (See, e.g., \cite{seneta1976}.) So without loss of
generality, for the case $\rho>0$, we may assume $\beta(\cdot) $ is
continuous and strictly increasing.

\item[2.] $\rho<0$: The transformation described in
\eqref{eqn:negRhoToPosRho} and \eqref{eqn:transformedCase}, show that
the pair $(X,Y)$ can be transformed to $(\widetilde X,Y)$ satisfying
$\rho>0$.

\item[3.] $\rho=0$: Suppose $\beta(\cdot) \in{\Pi_+(a)}$ after which we
consider $\beta\in\Pi_-(a)$. From \cite{dehaanresnick1979} as reviewed
in Section~\ref{subsubsec:fctClass}, there exists $\tilde\beta(t)$
which is continuous, strictly increasing and such that
$\beta-\tilde\beta=o(\alpha)$ so that the convergence of types theorem
allows us to replace $\beta$ by $\tilde\beta$. Assume this is done
which is tantamount to dropping the tilde. Then there are two cases to
consider.
\begin{longlist}[(a)]
\item[(a)]\label{item:1stSubcase} $\beta(\infty) =\infty.$

\item[(b)]\label{item:2ndSubcase} $\beta(\infty) <\infty.$
\end{longlist}
For~3(a) 
it is clear that $\beta(t)$ has the desired properties of being
continuous and strictly increasing to $\infty$. For~3(b), 
proceed as follows to transform $(X,Y)$: Define
%
\begin{eqnarray}\label{eqn:moreTransforms}
\widetilde X&=& \frac{1}{\beta(\infty) -X},\qquad
\tilde\beta(t)=\frac{1}{\beta(\infty) -\beta(t)},
\nonumber
\\[-8pt]
\\[-8pt]
\nonumber
\tilde\alpha(t)
&=& \frac{\alpha(t)}{(\beta(\infty) -\beta(t))^2}.
\end{eqnarray}
Then $\tilde\beta(t) \uparrow\infty$ is continuous and strictly
monotone and $\tilde\beta\in\Pi_+(\tilde\alpha)$ and
after some calculation we get
\begin{eqnarray*}
&& tP \biggl[ \frac{\widetilde X -\tilde\beta(t)}{\tilde\alpha(t)} \leq x,
\frac Yt>y \biggr]
\\
&&\qquad = tP \biggl[ \frac{X-\beta(t)}{\alpha(t)} \leq
\frac{x}{1+ \alpha(t)x/(\beta(\infty)-\beta(t))},\frac Yt>y \biggr]
\\
&&\qquad \to \mu\bigl( [-\infty,x]\times(y,\infty] \bigr)
\end{eqnarray*}
since $\tilde\beta\in\Pi_+(\tilde\alpha)$ implies $\tilde
\beta(t)/\tilde\alpha(t) \to\infty$ which is identical to
$(\beta(\infty)-\beta(t))/\alpha(t) \to\infty.$ Thus after the
transformation of $(X,Y)$ to $(\widetilde X,Y)$, case 3(b) is reduced
to case~3(a).

What if $\beta\in\Pi_-(a)$? Then define
\[
\widetilde X=-X,\qquad \tilde\beta(t)=-\beta(t),\qquad \tilde\alpha
(t)=\alpha(t),
\]
and $\tilde\beta\in\Pi_+(a)$ and this case reduces to the case
when $\beta\in\Pi_+(a)$ since
\begin{eqnarray*}
tP \biggl[ \frac{\widetilde X-\tilde\beta(t)}{\tilde\alpha(t)} \leq x,
\frac Yt>y \biggr] &=&  tP \biggl[ \frac{ X- \beta(t)}{ \alpha(t)}
\geq-x, \frac
Yt>y \biggr]
\\
&\hspace*{3pt} \to& \mu \bigl( [-x,\infty]\times(y,\infty] \bigr).
\end{eqnarray*}
\end{longlist}

\subsubsection{Summary}\label{subsubsec:summary}
When $\psi_2 \not\equiv0$, if we make the transformation
$X\mapsto\widetilde X$ and consider the analogue of
\eqref{eqn:basicConvStandardForm} for $(\widetilde X,Y) $, we can
standardize the $\widetilde X$-variable. If \mbox{$\psi_2 \equiv0$}, but
$\psi_1(c)=c^\rho,$ for $c>0$, $\rho\neq 0$, then for $\rho>0$,
$(\alpha^\leftarrow(X),Y)$ are a standardized pair and for $\rho<0$, $
( (1/\alpha)^\leftarrow(\widetilde X),Y )$ is a standardized pair.

When the limit $\mu$ is a product measure, $(\psi_1,\psi_2)\equiv
(1,0)$ and standardization is not possible; an example is given in
Section \ref{subsubsec:ParetoNoNo} and a proof of the assertion is easy
using the change of coordinate system techniques of
Section~\ref{sec:changeCoord}.

\subsection{Densities}\label{subsec:densities} In this section we see
what form the basic convergence takes when $(X,Y)$ has a density. Since
it is sufficient to suppose that the $Y$-variable has been transformed
to the standard case, for this section, we assume the following:
\begin{longlist}[3.]
\item[1.] The pair $(X,Y)$ has density $f(x,y)$.

\item[2.] The marginal density $f_Y(y)=\int_{-\infty}^\infty
f(x,y)\,dx$ of the $Y$-variable satisfies
\[
f_Y(y) = y^{-2},\qquad   y>1.
\]
Since we have densities, we assume the transformation to $Y$ being
standard renders $Y$ a Pareto random variable with unit shape
parameter.

\item[3.] The joint density $f(x,y)$ satisfies
%
\begin{equation}\label{eqn:basicConvDensities}
t^2 \alpha(t) f\bigl(\alpha(t)x+\beta(t),ty\bigr) \to g(x,y) \in L_1
\bigl( [-\infty,\infty]\times(0,\infty] \bigr),
\end{equation}
where the limit $g(x,y)\geq0$ is integrable, not identically zero
and satisfies for each fixed $v>0$,
%
\begin{equation}\label{eqn:requireDensity}
v^2g(u,v)\qquad\mbox{is a probability density in }u.
\end{equation}
\end{longlist}

\begin{Proposition}\label{prop:densityConvImpliesBasicConv}
With the assumptions just listed, \eqref{eqn:basicConvStandardForm}
holds with
\[
\mu\bigl( [-\infty,x]\times(y,\infty] \bigr)=\int_{u\leq
x}\int_{v>y}g(u,v)\, dv\,du,
\]
and $H(\infty)=\mu([-\infty,\infty]\times(1,\infty])=1.$
\end{Proposition}

\begin{pf}
We use standard notation for conditional densities. So for instance,
$f_{X|Y=v}(u|v)$ is the conditional density of $X$ given $Y=v$.

We need two facts:
\begin{longlist}[2.]
\item[1.] First we evaluate the integrand. For $v>0$,
\eqref{eqn:basicConvDensities} implies
%
\begin{equation}\label{eqn:conditDensityConv}
f_{(X-\beta(t))/\alpha(t)| Y/t =v } (u|v) \to v^2 g(u,v)\qquad (t
\to\infty).
\end{equation}
To see this, observe
\begin{eqnarray*}
f_{(X-\beta(t))/ \alpha(t)| Y/t =v } (u|v) &=&
\frac{f_{(X-\beta(t))/\alpha(t), Y/t } (u,v)} {f_{Y/t}(v)}
=\frac{t \alpha(t) f(\alpha(t)u +\beta(t),tv)}{tf_Y(tv)}
\\
&=& t^2 \alpha(t) v^2 f\bigl(\alpha(t)u +\beta(t),tv\bigr)\to v^2
g(u,v).
\end{eqnarray*}

\item[2.] We now show convergence of the integral. The function of $u$
\[
f_{(X-\beta(t))/\alpha(t)| Y/t =v } (u|v)
\]
is a probability density for fixed $v$.
\end{longlist}

Now write
\begin{eqnarray*}
&& tP \biggl[ \frac{X-\beta(t)}{\alpha(t)} \leq x, \frac Yt >y \biggr]
\\
&&\qquad =
t \int_{[v>y]} \biggl[ \int_{[u\leq x]} f_{(X-\beta(t))/\alpha(t)|  Y/t
=v } (u|v)\,du \biggr] f_{Y/t}(v)\,dv
\\
&&\qquad = \int_{[v>y]} \biggl[ \int_{[u\leq x]} f_{(X-\beta(t))/\alpha(t)|
Y/t =v } (u|v)\,du \biggr] v^{-2}\,dv.
\end{eqnarray*}
The integral inside the square bracket has an integrand which is a
family of probability densities in the variable $u$ (with $v$ fixed)
indexed by $t$ which converges to a limiting probability density
$v^2g(u,v).$ Hence by Scheff\'e's lemma (e.g., \cite{resnick1998book},
page~253)
\[
\biggl[ \int_{[u\leq x]} f_{(X-\beta (t))/\alpha(t)| Y/t =v } (u|v)\,du
\biggr] \to\int_{[u\leq x]}v^2 g(u,v)\,du.
\]
Now the square bracket term is a conditional
probability and hence is a function of $v$ bounded almost surely by 1.
So by dominated convergence, we have proven
\eqref{eqn:basicConvStandardForm} as required.

To check the last assertion that $H(\infty)=1$, note
\begin{eqnarray*}
\int_{-\infty}^\infty\int_{v>1} g(u,v)\,du \,dv&=&  \int_{v>1} v^{-2}
\biggl(\int_{-\infty}^\infty v^2 g(u,v)\,du \biggr)\,dv
\\
&=& \int_{v>1}
v^{-2}dv =1.
\end{eqnarray*} \upqed
\end{pf}

\mbox{Heffernan} and Tawn \cite{heffernantawn2004} assume that $(X,Y)$ have
been transformed to have Gumbel marginal distributions, that is,
$P(X\leq t)=P(Y\leq t)=\break \exp(-\exp(-t))$ for $t\in\mathbb{R}$ and that
for such $(X,Y)$
%
\begin{equation}\label{eqn:htConv}
tP \biggl[\frac{X-\tilde\beta(t)}{\tilde\alpha(t)} \leq x \Big| Y=t \biggr]
\end{equation}
converges to a nondegenerate limit distribution as $t\to\infty$, for
some scaling function $\tilde\alpha(\cdot)>0$ and centering function
$\tilde\beta(\cdot)\in\mathbb{R}$.

Thus we see that since~(\ref{eqn:conditDensityConv}) implies
\cite{heffernantawn2004} condition~(\ref{eqn:htConv}),
(\ref{eqn:basicConvDensities}) implies~(\ref{eqn:htConv}). This makes
explicit the link between our assumptions~(\ref{eqn:basicConv}) and
those of \mbox{Heffernan} and Tawn \cite{heffernantawn2004} under the above
conditions for densities. We have
\begin{eqnarray*}
P \biggl[\frac{X-\tilde\beta(t)}{\tilde\alpha(t)} \leq x \Big| Y=ty \biggr]
&=& \int_{u\leq x}f_{(X-\tilde\beta
(t))/\tilde\alpha(t) |Y/t =y }(u | y)\,du
\\
&\hspace*{3pt} \to&\int_{u\leq x}y^2 g(u,y)\,du,
\end{eqnarray*}
and letting $y=1$ gives
\[
P \biggl[\frac{X-\tilde\beta(t)}{\tilde\alpha(t)} \leq x \Big| Y=t \biggr]
\to\int_{u\leq x}g(u,1)\,du.
\]

\section{Characterizing the class of limit
measures}\label{sec:charLimits}
Assuming the $Y$-variable is standardized, what is the class of limits
in \eqref{eqn:basicConvStandardForm}?
We divide this issue in two parts, depending on whether the limit
measure $\mu$ is a product or not.

\subsection{The limit measure is a product}\label{subsec:prod} For
this case, there is not much discussion required since for any
distribution function $H(x)$ on $\mathbb{R}$, the limit
\[
\mu=H\times\nu_1\quad\mbox{or}\quad
\mu\bigl([-\infty,x]\times(y,\infty] \bigr)=H(x)y^{-1}
\]
is possible. To achieve this limit, suppose $X,Y$ are independent
random variables with $X$ having distribution $H$ and $Y$ being
standard Pareto. Then with $\beta(t)=0$ and $\alpha(t) =1$,
\eqref{eqn:basicConvStandardForm} is satisfied.

\subsection{The limit measure is not a product}\label{subsec:notprod}
When $\mu$ is not a product,
we change coordinate systems and transform $X$ to some $X^*$ and
assume $(X^*,Y)$ is a standard pair and
%
\begin{equation}\label{eqn:bothStandard}
tP \biggl[ \biggl( \frac{ X^* }{t} , \frac{Y}{t} \biggr) \in\cdot
\biggr] \stackrel{v}{\to} \mu_* (\cdot)\qquad\mbox{in } M_+ \bigl(
[0,\infty]\times(0,\infty] \bigr),
\end{equation}
where $\mu_*$ is a transformation of $\mu$ as described in Section
\ref{subsec:canWeStandardizeX}.

From \eqref{eqn:bothStandard}, we see that the distribution of $( X^*,
Y)$ is standard regularly varying with limit measure $\mu_*$ (see
\cite{basrakdavismikosch2000,resnickbook2006,resnick2004b}) on the
cone $[0,\infty]\times(0,\infty]$ and, therefore $\mu_*$ is homogeneous
of order -1:
\[
\mu_*(c\Lambda)=c^{-1}\mu_*(\Lambda)\qquad \forall c>0,
\]
where $\Lambda$ is a Borel subset of
$[0,\infty]\times(0,\infty]$. This means $\mu_* $ has a spectral
form. We pick a norm. Any norm would do but for convenience define
\[
\|(x,y)\|=|x|+|y|,\qquad (x,y)\in\mathbb{R}^2.
\]
Of course, when restricting attention to $[0,\infty]\times (0,\infty]$,
the absolute value bars can be dropped. Then the standard argument
using homogeneity (\cite{resnick1987}, Chapter~5), yields for $r>0$ and
$\Lambda$ a Borel subset of $[0,1)$,
%
\begin{eqnarray}\label{eqn:defS}
\nonumber && \mu_* \biggl\{ (x,y) \in[0,\infty]\times(0,\infty]\dvtx
x+y>r, \frac{x}{x+y} \in\Lambda \biggr\}
\\
&&\qquad =\mu_* \biggl\{ r(x,y)\in[0,\infty]\times(0,\infty]\dvtx
x+y>1, \frac{x}{x+y} \in\Lambda \biggr\}
\nonumber
\\[-8pt]
\\[-8pt]
\nonumber &&\qquad = r^{-1}\mu_* \biggl\{
(x,y)\in[0,\infty]\times(0,\infty]\dvtx x+y>1, \frac{x}{x+y} \in\Lambda
\biggr\}
\\
\nonumber &&\qquad =: r^{-1} S(\Lambda).
\end{eqnarray}
The Radon measure $S$ need not be a finite measure on $[0,1)$
but to guarantee that
%
\begin{equation}\label{eqn:defHStar}
H_*(x)=\mu_* \bigl( [0,x]\times (1,\infty] \bigr)
\end{equation}
is a probability measure, we need
%
\begin{equation}\label{eqn:condOnS}
\int_0^1(1-w)S(dw)=1.
\end{equation}
This will be clear from the following calculation to get the canonical
form of $H_*(x)$ for $x>0$:

Using \eqref{eqn:defS}, write for $x>0$,
%
\begin{eqnarray}\label{eqn:formForMuStar}
\nonumber && \mu_* \bigl([0,x]\times(y,\infty] \bigr)
\\
\nonumber &&\qquad  =\int\!\!\int_{\mathop{\mathop{0\leq rw \leq
x}\limits_{\mbox{\fontsize{8.36}{8}\selectfont{$r(1-w)>y$}}}}\limits_{\mbox{\fontsize{8.36}{8}\selectfont{$0\leq
w<1$}}}}
r^{-2}\,dr S(dw)
\\
&&\qquad = \int_{r=0}^\infty \biggl( \int_{\mathop{\mathop{0\leq w \leq
x/r}\limits_{\mbox{\fontsize{8.36}{8}\selectfont{$1-y/r>w$}}}}
\limits_{\mbox{\fontsize{8.36}{8}\selectfont{$0\leq w <1$}}}} S(dw) \biggr)r^{-2}\,dr
\\
\nonumber &&\qquad =\int_0^\infty S
\biggl(\biggl[0,\frac xr \wedge\biggl(1-\frac yr \biggr)\wedge
1\biggr) \biggr) r^{-2}\,dr
\\
\nonumber &&\qquad =\int_0^{{\infty}} S \bigl(\bigl[0,xv
\wedge(1-yv)\wedge1\bigr) \bigr)\,dv.
\end{eqnarray}

Integrating the double integral in reverse order yields the
alternate expression
%
\begin{eqnarray}\label{eqn:altForm}
\nonumber && \mu_* \bigl([0,x]\times(y,\infty] \bigr)
\\
&&\qquad = \int_{ w
\in[0,1)} \biggl( \int_{ y/(1-w)<r\leq  x/w } r^{-2}\,dr \biggr) S(dw)
\nonumber
\\[-8pt]
\\[-8pt]
\nonumber
&&\qquad = \int_{w\in[0,1)} \bigl(
(1-w)y^{-1} -wx^{-1} \bigr)_+ S(dw)
\\
\nonumber &&\qquad = y^{-1}\int_0^{x/(x+y)} (1-w)S(dw)
-x^{-1}\int_0^{x/(x+y)}wS(dw).
\end{eqnarray}

Conclusion: The class of limits $\mu_*$ or conditional limits
\[
H_*(x)=\lim_{t \to\infty}P \biggl[\frac{X^*}{t} \leq x\Big|Y>t\biggr]
\]
is indexed by Radon measures $S$ on $[0,1)$ satisfying the
integrability condition~\eqref{eqn:condOnS}.

\begin{Example*}
As an example, suppose $S$ is uniform on $[0,1)$: $S(dw)=\frac{dw}{c},$
where $c$ is chosen so that \eqref{eqn:condOnS} is satisfied: $\int_0^1
\frac wc \,dw {=1}$ which implies $c=1/2$. This yields
\[
\mu_* \bigl([0,x]\times(y,\infty] \bigr) =\frac{x}{x+y} \biggl[ \frac2y
- \biggl( \frac{1+ x/y }{x+y} \biggr) \biggr]
\]
and setting $y=1$ we get a Pareto distribution
\[
H_*(x)=\frac{x}{1+x}=1-\frac{1}{1+x},\qquad  x>0.
\]
\end{Example*}

\section{Random norming}\label{sec:randomNorming}
In \cite{heffernantawn2004}, it was necessary to normalize $X$ by a
function of the precise value of $Y$ occurring with $X$ to achieve
nondegeneracy of the limiting conditional distribution. Motivated by
this, we consider how to normalize the \mbox{$X$-}variable with a function of
$Y$ rather than a deterministic affine transformation, {using functions
of the threshold $t$ in~(\ref{eqn:conditVersion})}. This leads to a
product form limit in all cases.

It is significant that normalizing by using functions of the threshold
$t$ in~(\ref{eqn:conditVersion}) does not result in a product limit in
all cases, but that the inclusion of the precise value of $Y$ occurring
with $X$ adds enough detail to the normalization to allow the limit
always to factorize. In statistical applications the factorization of
the limit distribution will constitute a welcome simplification of
models based on this limiting form. Indeed, the statistical model of
\mbox{Heffernan} and Tawn \cite{heffernantawn2004}relies on such factorization
to ensure that the residuals formed by normalizing observed values of
$X$ by functions of the observed values of $Y$ are independent of the
$Y$ values.

We discuss this random normalization in two stages:
\begin{itemize}
\item The $X$-variable can be standardized and the limit
in~\eqref{eqn:basicConvStandardForm} is not a product.

\item The limit measure $\mu$ in~\eqref{eqn:basicConvStandardForm} is a
product measure.
\end{itemize}

\subsection{The $X$-variable can be standardized and the limit measure $\mu$ is not a
product}\label{subsec:randomNormingXStandard}
We suppose $X$ can be transformed to $X^*$ so that $(X^*,Y)$ is a
standardized pair and \eqref{eqn:bothStandard} holds with limit
measure $\mu_*$.
As in Section
\ref{subsec:notprod}, let $S$ be the spectral measure of $\mu_*$. Then
we have the following result which forms the basis of the estimation
procedure proposed in \cite{heffernantawn2004}.

\begin{Proposition}\label{prop:randomNorming}
If \eqref{eqn:bothStandard} holds, then
%
\begin{equation}\label{eqn:randomNorming}
tP \biggl[ \biggl( \frac{X^*}{Y}, \frac Yt \biggr) \in\cdot \biggr]
\stackrel{v}{\to} G\times\nu_1\qquad\mbox{in }M_+
\bigl([0,\infty]\times(0,\infty] \bigr),
\end{equation}
where for $x>0$
%
\begin{equation}\label{eqn:addin}
\nu_1 ((x,\infty] ) =x^{-1}\quad\mbox{and}\quad G(x)=\int_0^{x/(1+x)}
(1-w)S(dw).
\end{equation}
This means
\[
P \biggl[ \frac{X^*}{Y} \leq x\Big| Y>t \biggr] \to G(x), \qquad  x>0.
\]

Conversely, if \eqref{eqn:randomNorming} holds, then so does
\eqref{eqn:bothStandard}.
\end{Proposition}

\begin{pf}
This proof is discussed in Theorem 2.1 of
\cite{maulikresnickrootzen2002}. The outline of the argument is as
follows. Applying the map $T_1(x,y)= (\frac xy,y )$ to
\eqref{eqn:bothStandard} yields after a compactification argument that
\[
tP \biggl[ \biggl( \frac{X^*}{Y}, \frac Yt \biggr) \in\cdot \biggr]
\stackrel{v}{\to} \mu_* \circ T_1^{-1}.
\]
So the limit
evaluated on $[0,x]\times(y,\infty]$ is
\begin{eqnarray*}
&& \mu_* \biggl\{(u,v)\dvtx \frac uv \leq x, v>y\biggr\}
\\
&&\qquad =
y^{-1}\mu_* \biggl\{(u,v)\dvtx \frac uv \leq x, v>1\biggr\}
\\
&&\qquad = y^{-1}\int\!\!\int_{\mathop{rw/(r(1-w)) \leq
x}\limits_{\mbox{\fontsize{8.36}{8}\selectfont{$r(1-w)>1$}}}}
r^{-2}\,dr\,S(dw)
\\
&&\qquad =y^{-1} \int_{w \leq x/(1+x)} \biggl( \int_{r>1/(1-w)}
r^{-2}\,dr\biggr)S(dw)
\\
&&\qquad =y^{-1} \int_0^{x/(1+x)} (1-w)S(dw).
\end{eqnarray*}

The converse proceeds similarly using the map
$T_2(x,y)=(xy,y)=T_1^{-1}(x,y).$
\end{pf}

\subsection{The limit measure $\mu$ is a product measure}\label{subsec:randomNormingLimitProduct}
Now we suppose \eqref{eqn:basicConvStandardForm} holds
with $\mu=H\times\nu_1$. In this case, from Proposition
\ref{prop:whenMuProd}, \eqref{eqn:varyAlpha} and \eqref{eqn:varyBeta}
hold with \mbox{$\psi_1(x)\equiv1$,} $\psi_2(x)\equiv0.$

\begin{Proposition}\label{prop:whoKnows}
If,
%
\begin{equation}\label{eqn:reviewStandard}
tP \biggl[ \frac{X-\beta(t)}{\alpha(t)} \leq x, \frac Yt>y \biggr] \to
H(x)y^{-1}\qquad (x\in\mathbb{R}, y>0)
\end{equation}
for a nondegenerate probability distribution function $H(x)$, then also
%
\begin{equation}\label{eqn:tastySandwich}
tP \biggl[ \frac{X-\beta(Y)}{\alpha(Y)} \leq x, \frac Yt>y \biggr] \to
H(x)y^{-1}\qquad (x\in\mathbb{R}, y>0 )
\end{equation}
and
\[
P \biggl[ \frac{X-\beta(Y)}{\alpha(Y)} \leq x\Big| Y>t \biggr] \to H(x).
\]

Conversely, if \eqref{eqn:tastySandwich} holds and $\alpha(\cdot)$ and
$\beta(\cdot)$ satisfy \eqref{eqn:varyAlpha}, \eqref{eqn:varyBeta}
locally uniformly with $\psi_1(x)\equiv1,$ and $\psi_2(x)\equiv0,$ then
\eqref{eqn:reviewStandard} also holds.
\end{Proposition}

\begin{pf}
For any $K>y>0$ we have
\begin{eqnarray*}
&& tP \biggl[ \frac{X-\beta(Y)}{\alpha(Y)} \leq x, \frac Yt \in(y,K]
\biggr]
\\
&&\qquad =tP \biggl[ \frac{X-\beta(t)}{\alpha(t)} \leq \frac{\alpha(t
 Y/t) }{\alpha(t)}x +\frac{\beta(t  Y/t) -\beta(t)
}{\alpha(t)}, \frac Yt \in(y,K]
\biggr]
\end{eqnarray*}
and because of local uniform convergence in \eqref{eqn:varyAlpha} and
\eqref{eqn:varyBeta}, this converges to
\[
\mu\bigl( [-\infty,x]\times (y,K] \bigr)=H(x) (y^{-1}-K^{-1} ).
\]
Therefore
\begin{eqnarray*}
\liminf_{t \to\infty} tP \biggl[ \frac{X-\beta(Y)}{\alpha(Y)} \leq x,
\frac Yt >y \biggr] &\geq& \liminf_{t \to\infty} tP \biggl[
\frac{X-\beta(Y)}{\alpha(Y)} \leq x, \frac Yt \in(y,K]
\biggr]
\\
&=& H(x) (y^{-1}-K^{-1} ).
\end{eqnarray*}
Since this is true for all $K>y$, we have
\[
\liminf_{t \to\infty} tP \biggl[ \frac{X-\beta(Y)}{\alpha(Y)} \leq x,
\frac Yt >y \biggr] \geq H(x)y^{-1}.
\]
Also,
\begin{eqnarray*}
\limsup_{t \to\infty} tP \biggl[ \frac{X-\beta(Y)}{\alpha(Y)} \leq
x, \frac Yt >y \biggr]
&\leq& \lim_{t \to\infty} tP \biggl[
\frac{X-\beta(Y)}{\alpha(Y)} \leq x, \frac Yt \in(y,K]
\biggr]
\\
&&{}  + \limsup_{t \to\infty} tP \biggl[ \frac Yt >K
\biggr]
\\
&=& H(x) ({y^{-1}-K^{-1}} ) + K^{-1}.
\end{eqnarray*}
Letting $K \to\infty$ provides the other half of the sandwich and
\eqref{eqn:tastySandwich} is proven.

For the converse, write
\begin{eqnarray*}
&& tP \biggl[\frac{X-\beta(t)}{\alpha(t)} \leq x, \frac Yt \in
(y,K] \biggr]
\\
&&\qquad = tP \biggl[ \frac{X-\beta(Y)}{\alpha(Y)} \leq
\frac{\alpha(t)}{\alpha(Y)}x + \frac{\beta(t)-\beta(Y)}{\alpha(Y)}
,\frac Yt \in(y,K]  \biggr].
\end{eqnarray*}
Proceed as before using uniform convergence.
\end{pf}

\section{Connection to multivariate extreme value theory and asymptotic independence}
\label{sect:ConnectionToAsympIndpce}
We now make some comments on the relationship between our conditioned
limit condition \eqref{eqn:basicConvStandardForm} and multivariate
extreme value theory.

Suppose the distribution of $(X,Y)$ is in the domain of attraction of
a multivariate extreme value distribution. This means that for i.i.d.
replicates $\{(X_i,Y_i), i \geq1\}$ of $(X,Y)$ there exist centering
$b_j(t) \in\mathbb{R}$ and scaling $a_j(t)>0$ functions, $j=1,2$, and
%
\begin{equation}\label{doa}
P\biggl[\frac{\bigvee_{i=1}^n X_i -b_1(n)}{a_1(n)} \leq x,
\frac{\bigvee_{i=1}^n Y_i -b_2(n)}{a_2(n)} \leq y\biggr] \to G(x,y),
\end{equation}
where $G$ is a multivariate extreme value distribution. Let the
marginal distributions of $G$ be $G_j$, $j=1,2$. \textit{Asymptotic
independence} means $G(x,y)=G_1(x)G_2(y)$.

Define
\begin{eqnarray*}
U_1(x)&=& \frac{1}{P[X>x]},\qquad  U_2(y)=\frac{1}{P[Y>y]},
\\
\chi_j(x)&=&  \biggl( \frac{1}{-\log G_j} \biggr)^\leftarrow(x),\qquad
x>0, j=1,2,
\\
G^*(x,y)&=& G(\chi_1(x),\chi_2(y)),\qquad  x>0,y>0.
\end{eqnarray*}
According to Resnick \cite{resnick1987}, Proposition~5.10, page~265, we
can standardize the condition \eqref{doa} by transforming $(X,Y)\mapsto
(X^*,Y^*)=(U_1(X), U_2(Y))$ and then
%
\begin{equation}\label{doa:standard}
P\biggl[\frac{\bigvee_{i=1}^n X^*_i}{n} \leq x, \frac{\bigvee_{i=1}^n
Y^*_i }{n} \leq y\biggr] \to G^*(x,y),
\end{equation}
and $G^*$ is max-stable. From \cite{resnick1987}, Proposition 5.15,
page 277 and \cite{resnickbook2006}, Section 6.1, this is equivalent to
marginal convergence and multivariate regular variation of the
distribution of $(X^*,Y^*)$:
%
\begin{equation}\label{multRegVar}
tP\biggl[ \biggl( \frac{X^*}{t}, \frac{Y^*} {t} \biggr) \in\cdot\biggr]
\stackrel{v}{\to} \nu^* (\cdot),
\end{equation}
in $M_+([0,\infty]^2 \setminus\{\bzero\})$. Here
$\nu^* $ is a Radon measure on $[0,\infty]^2 \setminus\{\bzero\}$
satisfying
%
\begin{equation}\label{scaling}
\nu^*(t \cdot)= t^{-1}\nu^*(\cdot).
\end{equation}
Asymptotic independence means
\begin{eqnarray*}
\nu^*([0,x]\times[0,y])^c&=& -\log G^*(x,y)
=-\log G^*(x,\infty)-\log G^*(\infty,y)
\\
&=& \nu^*\bigl((x,\infty]\times[0,\infty]\bigr)
+\nu^*\bigl([0,\infty]\times(y,\infty] \bigr),
\end{eqnarray*}
and $\nu^*$ concentrates on the lines $\{ (x,0)\dvtx x > 0
\}\cup\{(0,y)\dvtx y > 0\}.$

Suppose the domain of attraction condition \eqref{multRegVar} holds but
asymptotic independence does \textit{not} hold. Condition~\eqref{multRegVar} implies for $x>0, y>0$,
\[
tP\biggl[\frac{X^*}{t} \leq x,\frac{Y^*}{t} >y\biggr] \to\nu^* \bigl(
[0,x]\times(y,\infty] \bigr)
\]
and we claim for
fixed $y>0$, $ \nu^* ( [0,x]\times(y,\infty] )$ is not
degenerate in $x$. This follows, for instance, from~\eqref{scaling}.
Conclusion: the domain of attraction condition \eqref{multRegVar} in
standard form without asymptotic independence implies that $(X^*,Y^*)$
satisfy~\eqref{eqn:basicConvStandardForm}. Condition
\eqref{eqn:basicConvStandardForm} is equivalent to vague convergence
on the cone $[0,\infty]\times(0,\infty]$ while the regular variation
condition \eqref{multRegVar} gives vague convergence on the bigger
cone $[0,\infty]^2\setminus\{\bzero\}$.

Suppose \eqref{multRegVar} holds \textit{with} asymptotic independence.
Consider \eqref{eqn:basicConvStandardForm} with $X^*/t$ in place of
$(X-\beta(t))/\alpha(t)$. The nondegeneracy condition in
\eqref{eqn:basicConvStandardForm} fails because for fixed $y>0$,
$\mu([-\infty,x]\times (y,\infty] )=\nu^* ([-\infty,x]\times(y,\infty]
) $ concentrates all mass at $x=0$. If one wants
\eqref{eqn:basicConvStandardForm} to hold, one must make an additional
assumption beyond the domain of attraction condition \eqref{multRegVar}
and the $X^*$ variable in \eqref{multRegVar} must be normalized
differently. For a simple particular case which is somewhat familiar,
consider the following: Suppose we assume the condition
\eqref{multRegVar} with asymptotic independence and in addition we
assume that $X^*$ can be normalized by $\alpha(t)$ instead of by $t$,
so that \eqref{eqn:basicConvStandardForm} holds in the form
%
\begin{equation}\label{basicParticular}
tP\biggl[\frac{X^*}{\alpha(t)} \leq x,\frac{Y^*}{t} >y\biggr]
\to\mu\bigl([0,x]\times(y,\infty]\bigr),\qquad  x>0, y>0.
\end{equation}
From \eqref{basicParticular} and \eqref{multRegVar}, we have for $0<a<b
\leq\infty$ and $y>0$
\begin{eqnarray*}
tP\biggl[\frac{X^*}{\alpha(t)} \in(a,b],\frac{Y^*}{t} >y\biggr]
&\to& \mu\bigl((a,b]\times(y,\infty]\bigr),
\\
tP\biggl[\frac{X^*}{t} \in(a,b],\frac{Y^*}{t} >y\biggr] &\to& 0.
\end{eqnarray*}

We claim that $t/\alpha(t) \to\infty$ so that $\alpha(\cdot) $ is of
smaller order than $t$. If not, there exist $t_n \to\infty$ and $0\leq
c <\infty$ and $t_n/\alpha(t_n) \to c$. From the nondegeneracy
condition in \eqref{eqn:basicConvStandardForm}, we may pick $0<a<b$
such that $\mu((a,b]\times(1,\infty])>0$. Then
\[
0<\mu\bigl((a,b]\times(1,\infty]\bigr) =\lim_{n \to\infty}
t_nP\biggl[\frac{X^*}{t_n} \in
\biggl(\frac{\alpha(t_n)}{t_n}a,\frac{\alpha(t_n)}{t_n}b\biggr],\frac{Y^*}{t_n}
>1\biggr]=0
\]
giving a contradiction. So $\alpha(\cdot)$ is of smaller order than $t$
and we have the situation of \textit{hidden regular variation}
\cite{heffernanresnick2004,maulikresnick2003b,resnick2002a}; that is,
the regular variation condition \eqref{multRegVar} holds on the big
cone $[0,\infty]^2\setminus \{\bzero\}$ but a different regular
variation condition holds on the smaller cone
$[0,\infty]\times(0,\infty]$.

To summarize: The multivariate extreme value paradigm without asymptotic
independence subsumes our conditioned limit condition
\eqref{eqn:basicConv}. However, in the presence of asymptotic
independence, the multivariate extreme value condition is refined by
\eqref{eqn:basicConv} which uses a more delicate normalization to
track mass into the part of the distributional tail where the
conditioning variable $Y$ is large.

\section{Examples}\label{sec:egs}
We give examples to illustrate some intricacies.

\subsection{Bivariate normal}\label{subsec:bivarNormal}
Suppose $N_1,N_2$ are i.i.d. $N(0,1)$ random variables and $|\rho|\leq
1$. Define $(X,Y)=(\sqrt{1-\rho^2}N_1+\rho N_2,N_2)$ which is a
bivariate normal vector with means 0, variances 1 and correlation
$\rho$. Denote the standard normal distribution function by $N(x)$.
Recall (e.g., from \cite{resnick1987}, page~71) that we may set
%
\begin{eqnarray}\label{eqn:normingConstants}
a(t)&=& \frac{1}{\sqrt{2\log t}}, \nonumber
\\[-8pt]
\\[-8pt]
\nonumber  b(t)&=&  \biggl(\frac{1}{1-N} \biggr)^\leftarrow(t)=
\sqrt{2\log t} -\frac{(1/2) (\log \log t+\log4\pi)}{\sqrt{2\log t}}+
o(a(t)),
\end{eqnarray}
and then for $ x\in\mathbb{R}$,
\[
\lim_{t\to\infty} tP\biggl[\frac{N_1-b(t)}{a(t)}>x\biggr] = e^{-x}.
\]

\subsubsection{Conditional limits for $(X,Y)$}\label{subsubsec:xy}
We begin by discussing the following result learned from
\cite{abdousfougeresghoudi2005}. Suppose $N(x)$ is the standard
normal distribution function and $n(y)$ is its density. Then
%
\begin{equation}\label{eqn:xy}
tP\biggl[X-\rho b(t) \leq x, \frac{Y-b(t)}{a(t)} >y\biggr]\to
N\bigl(x/\sqrt{1-\rho^2}\bigr) e^{-y},
\end{equation}
or standardizing the $Y$-variable,
%
\begin{equation}\label{eqn:xyStandard}
tP\biggl[X-\rho b(t) \leq x, \frac{b^\leftarrow(Y)}{t} >y\biggr]\to
N\bigl(x/\sqrt{1-\rho^2}\bigr) y^{-1}.
\end{equation}
Here we claimed $\beta(t) =\rho b(t)$ and $ \alpha(t) =1.$ It is well
known (e.g., \cite{resnick1987}, page~71) that $b(\cdot)
\in\Pi(a(\cdot))$ and therefore
%
\begin{eqnarray}\label{eqn:goesToZero}
\frac{\beta(tc)-\beta(t)}{\alpha(t)} &=& \rho \bigl(b(tc)-b(t) \bigr)
\nonumber
\\[-8pt]
\\[-8pt]
\nonumber
&=&  \rho\frac{ (b(tc)-b(t) )}{a(t)}a(t) \sim \rho\log c \cdot a(t) \to0.
\end{eqnarray}
Thus $\psi_2(x)$ in \eqref{eqn:varyBeta} is identically 0 and
$\psi_1(x) \equiv1.$

We now see why \eqref{eqn:xy} and
\eqref{eqn:xyStandard} are true. We write,
\begin{eqnarray*}
&& tP\biggl[X-\rho b(t) \leq x, \frac{Y-b(t)}{a(t)}>y\biggr]
\\
&&\qquad = tP\biggl[\sqrt{1-\rho^2} N_1+\rho N_2 -\rho b(t)\leq x,
\frac{N_2-b(t)}{a(t)}>y\biggr]
\\
&&\qquad =\int_{a(t)y+b(t)}^\infty P\bigl[\sqrt{1-\rho^2} N_1+\rho{s}
-\rho b(t)\leq
x\bigr]tn(s)\,ds
\\
&&\qquad =\int_{y}^\infty P\bigl[\sqrt{1-\rho^2} N_1+\rho
\bigl(a(t)u+b(t) \bigr) -\rho b(t)\leq x\bigr]
\\
&&\phantom{\qquad =\int_{y}^\infty} {} \times
ta(t)n\bigl(a(t)u+b(t)\bigr)\,du
\\
&&\qquad \sim \int_y^\infty P\bigl[\sqrt{1-\rho^2}N_1 \leq x-{\rho}a(t)u\bigr]e^{-u}\,du
\end{eqnarray*}
since $ta(t)n(a(t)u+b(t)) \to e^{-u}.$ Using the fact that $a(t) \to0$,
we get convergence to
\[
\to \int_y^\infty P\bigl[\sqrt{1-\rho^2}N_1 \leq x\bigr]e^{-u}\, du
=N\bigl( x/\sqrt{1-\rho^2}\,\bigr)e^{-y},
\]
as claimed.

Conclusion: The limit measure is a product measure, $(\psi_1,\psi_2)
\equiv(1,0)$ and $\alpha(t)=1$. We have an illustration of
Proposition~\ref{prop:whenMuProd}.

\subsubsection{Exponential marginals for~$X$}\label{subsubsec:marginalsx}
In light of the standard form result \eqref{eqn:xyStandard} it is
tempting to look at limits for $(b^\leftarrow(X),b^\leftarrow(Y))$
but this turns out not to work. The reason for this is explored in
Section~\ref{subsubsec:ParetoNoNo}. Instead, following
\cite{heffernantawn2004}, we consider $(\log b^\leftarrow(X), \log
b^\leftarrow(Y))$. Thus we can transform $X$ to have exponential
marginals but not Pareto marginals.

We show the standard form
%
\begin{eqnarray}\label{eqn:expMarginals}
&& tP \biggl[ \frac{\log b^\leftarrow(X) -\log b^\leftarrow(\rho
b(t))}{\rho b(t)} \leq x, \frac{b^\leftarrow(Y)}{t}>y \biggr]
\nonumber
\\[-8pt]
\\[-8pt]
\nonumber
&&\qquad  \to N
\biggl(\frac{x}{\sqrt{\rule{0pt}{10pt}\smash{1-\rho^2}}} \biggr)y^{-1}.
\end{eqnarray}

The verification of
\eqref{eqn:expMarginals} needs the following lemma.
\begin{Lemma}\label{lemma:Pi}
The function
\[
V(t){:=}-\log\widebar N(\log t)=\log b^\leftarrow(\log t) \in\Pi(\log
t)
\]
is $\Pi$-varying with auxiliary function
$g(t)=\log t.$
\end{Lemma}

\begin{pf}
To prove membership in the $\Pi$-class, it suffices according to de
Haan \cite{dehaan1976} (see alternatively \cite{resnick1987}, page~30),
to show $V'(t) \in RV_{-1}$ and then the auxiliary function can be
taken to be $tV'(t)$. So it suffices to show
\[
(-\log\widebar N (\log t) )' \sim\frac{\log t}{t} \in RV_{-1}.
\]
The derivative is
\[
\frac{n(\log t)t^{-1}}{\widebar N(\log t)}\sim \frac{n(\log
t)t^{-1}}{n(\log t)/\log t}=t^{-1}\log t \in RV_{-1}.
\]\upqed
\end{pf}

To show \eqref{eqn:expMarginals}, we use \eqref{eqn:xyStandard} and
the Delta method. The left-hand side of \eqref{eqn:expMarginals} is
\begin{eqnarray*}
&& tP \biggl[ \frac{V ( e^{X-\rho b(t)}e^{\rho b(t)} ) -V ( e^{\rho
b(t)} )} { g (e^{\rho b(t)} )}\leq x, \frac{ b^\leftarrow(Y)}{t}>y
\biggr]
\\
&&\qquad \to P\bigl[\log e^{{N_1} \sqrt{\rule{0pt}{8pt}\smash{1-\rho^2}}}\leq x\bigr]y^{-1} =N
\biggl(\frac{x}{\sqrt{\rule{0pt}{10pt}\smash{1-\rho^2}}} \biggr) y^{-1}.
\end{eqnarray*}

Here is the conditional form of \eqref{eqn:expMarginals}, where $X$ is
transformed to have exponential marginals:
\begin{eqnarray*}
&& \lim_{t \to\infty} P \biggl[ \frac{\log b^\leftarrow(X) -\log
b^\leftarrow(\rho
b(t))}{\rho b(t)} \leq x\Big|Y>b(t) \biggr]
\\
&&\qquad =\lim_{t \to\infty} P \biggl[ \frac{\log b^\leftarrow(X) -\log
b^\leftarrow(\rho t)}{\rho t} \leq x\Big|Y>t \biggr] =N
\biggl(\frac{x}{\sqrt{\rule{0pt}{10pt}\smash{1-\rho^2}}} \biggr).
\end{eqnarray*}
The conditional form of \eqref{eqn:xyStandard}, where the marginal
distribution is normal, has the same limit:
\[
\lim_{t \to\infty}P[X-\rho b(t) \leq x|Y>b(t)] = \lim_{t \to
\infty}P[X-\rho t \leq x|Y>t]=N \biggl(\frac{x}{\sqrt{\rule{0pt}{10pt}\smash{1-\rho^2}}}
\biggr).
\]
This result seems natural when one observes that the normal
distribution is in the domain of attraction of the Gumbel distribution.

After transformation of $X$ to exponential marginals, we
have for \eqref{eqn:expMarginals}
\[
\beta(t) =-\log\widebar N(\rho b(t)),\qquad \alpha(t) =\rho b(t),
\]
and again $\psi_2(t)=0,$ since
\begin{eqnarray*}
\frac{\beta(tc)-\beta(t)}{\rho b(t)}&=& \frac{ \log (  \widebar N(\rho
b(tc))/\widebar N(\rho b(t)) ) }{\rho b(t)}\sim \frac{ \log (  n(\rho
b(tc)) / n(\rho b(t)) ) }{\rho b(t)}
\\
&\sim& \frac{\log e^{(\rho^2/2) (b^2(tc)-b^2(t) )}}{\rho b(t)}
=\frac{\rho^2}{2} \bigl( b(tc)-b(t) \bigr)\frac{ (b(tc)+b(t)
)}{\rho b(t)}
\\
&\sim& \rho\bigl(b(tc)-b(t) \bigr)\to0,
\end{eqnarray*}
using the same argument as in \eqref{eqn:goesToZero}. (This provides
another illustration of Proposition~\ref{prop:whenMuProd}.)

\subsubsection{Why $X$ cannot be transformed to
Pareto}\label{subsubsec:ParetoNoNo}
It is noteworthy that one cannot transform $X$ to have Pareto
marginals and expect the analogue of
\eqref{eqn:xy} to hold. Here is the explanation which
also relates to the discussion in Section
\ref{sec:changeCoord}.

Suppose for some
choice of centering and scaling $\alpha_2(t)>0, \beta_2 (t)\in
\mathbb{R} $ we have
%
\begin{equation}\label{eqn:BetterNotConv}
\lim_{t \to\infty} tP \biggl[ \frac{b^\leftarrow(X)-\beta_2(t)}
{\alpha_2(t)} \leq x, \frac{b^\leftarrow(Y)}{t}>y \biggr]
\end{equation}
exists and is nondegenerate in the sense of condition (iii) stated at
the beginning of Section~\ref{sec:basicDefResults}. This expression
\eqref{eqn:BetterNotConv} equals
%
\begin{equation}\label{eqn:leadsToNoNo}
\lim_{t \to\infty} P\biggl[X-\rho b(t) \leq b
\bigl(\alpha_2(t)x+\beta_2(t) \bigr) -\rho b(t),
\frac{b^\leftarrow(Y)}{t}>y\biggr]
\end{equation}
and from \eqref{eqn:xy} we would have for some nondecreasing limit
$\psi(x)$, that as $t \to\infty$,
%
\begin{equation}\label{eqn:interiorConv}
b \bigl(\alpha_2(t)x+\beta_2(t) \bigr) -\rho b(t) \to\psi(x).
\end{equation}
Furthermore, the limit in \eqref{eqn:BetterNotConv} would have to be
%
\begin{equation}\label{eqn:formLimit}
N \biggl(\frac{\psi(x)}{\sqrt{\rule{0pt}{10pt}\smash{1-\rho^2}}} \biggr)y^{-1}.
\end{equation}
Inverting \eqref{eqn:interiorConv}, we would need
\[
\frac{b^\leftarrow(y+\rho b(t)) -\beta_2(t)}
{\alpha_2(t)} \to\psi^\leftarrow(y).
\]
Changing variables leads to
\[
\frac{b^\leftarrow(\log tx)) -\beta_2(b^\leftarrow(\log t/\rho))}
{\alpha_2(b^\leftarrow(\log t/\rho))} \to\psi^\leftarrow (\log x).
\]
If $\psi^\leftarrow$ is not constant, then (\cite{gelukdehaan1987},
page~16)
\[
b^\leftarrow\circ\log= \biggl(\frac{1}{1-N} \biggr) \circ\log
\]
is either regularly varying with positive index or it is $\Pi$-varying.
Neither of these possibilities is true. If $\psi^\leftarrow$ is
constant, then the limit \eqref{eqn:formLimit} fails the nondegeneracy
assumptions.

So assuming the nondegenerate limit exists in \eqref{eqn:BetterNotConv}
leads to a contradiction. This illustrates the restrictions in our
ability to standardize the $X$ variable discussed in Section
\ref{subsec:canWeStandardizeX}.

\subsection{Heavy tailed examples}\label{subsec:HT}
{In this section, we present examples of heavy tailed random
variables possessing asymptotic independence.

\subsubsection{Mixture of independent standard regularly varying
random variables \textup{I:} positive $\rho$}\label{subsubsec:srvPower}
Suppose nonnegative random variables $(U,V)$ have a joint distribution
which is standard regularly varying; that is, there is a limit measure
$\nu$ on $[0,\infty]^2 \setminus\{\bzero\}$ such that
\[
tP \biggl[ \biggl( \frac Ut, \frac Vt \biggr)\in\cdot \biggr]
\stackrel{v}{\to} \nu
\]
in $M_+ ([0,\infty]^2 \setminus\{\bzero\} )$. For example, $(U,V)$
could be max-stable (\cite{resnick1987}, Chapter~5),
\cite{dehaanferreira2006} with exponent $\nu$. Suppose $(U_i,V_i),
i=1,2$, are i.i.d. copies of $(U,V)$. For $0<p<1$, define
%
\begin{equation}\label{eqn:defXY}
(X,Y)=B(U_1,V_1^p)+(1-B)(U_2^p,V_2),
\end{equation}
where
$P[B=0]=P[B=1]=\frac12,$ and $B$ is independent of $(U_i,V_i),
i=1,2.$

Observe that for any $x>0, y>0$
%
\begin{eqnarray}\label{eqn:asyIndep}
\nonumber && tP \biggl\{\biggl[\frac Xt \leq x, \frac Yt \leq y\biggr]^c
\biggr\}
\\
&&\qquad =\frac t2 P \biggl[\frac{U_1}{t}>x \mbox{ or }
\frac{V^p_1}{t}>y \biggr] + \frac t2 P \biggl[\frac{U^p_2}{t}>x \mbox{
or }\frac{V_2}{t}>y \biggr]
\\
\nonumber &&\qquad =\frac t2 P[U_1>tx]+o(1)+\frac t2 P[V_2>ty]+o(1) \to
\frac12 (x^{-1}+y^{-1}).
\end{eqnarray}
So $(X,Y)$ is standard regularly varying, in a domain of attraction of
a multivariate extreme value distribution, and possesses
asymptotic independence. The asymptotic independence holds
even if $(U,V)$ has no
asymptotic independence.

Now observe that
%
\begin{eqnarray}\label{eqn:defnu}
\nonumber && tP\biggl[ \frac{X}{t^p} \leq x, \frac Yt >y\biggr]
\\
&&\qquad =  \frac
t2 P[U_1 \leq t^p x, V_1^p>ty]+\frac t2 P[U_2^p \leq t^p
x,V_2>ty]
\nonumber
\\[-8pt]
\\[-8pt]
\nonumber
&&\qquad =  \frac t2 P[U_1 \leq t^p x, V_1>t^{1/p}y{^{1/p}}] +\frac
t2P[U_2\leq tx^{1/p} , V_2>ty]
\\
\nonumber &&\qquad \to 0 + \frac12 \nu \bigl( [0,x^{1/p}]\times(y,\infty]
\bigr) =:\mu\bigl([0,x]\times(y,\infty]\bigr).
\end{eqnarray}

If $(U,V)$ possess asymptotic independence, then $\nu((0,\infty]^2 )=0$
and the nondegeneracy assumption for $\mu$ stated in
\eqref{eqn:basicConvStandardForm} fails since for fixed $y>0$, the
function of $x$ given by $\nu ( [0,x^{1/p}]\times(y,\infty] )$
concentrates at $x=0$. So for this example, $(X,Y)$ is standard
regularly varying, asymptotically independent and provided $(U,V)$ does
not possess asymptotic independence, we can refine the asymptotic
independence to get the limit in \eqref{eqn:basicConvStandardForm}.
This gives an example of case (i) of \eqref{eqn:psi2} with $\rho=p$,
$\beta(t)=(1/\rho) \alpha(t)=t^p.$ The conditional limit distribution
can most simply be written as
\[
\lim_{t \to\infty}P\biggl[\frac{X}{t^p} \leq x\Big|Y>t\biggr]= \frac12
\nu\bigl( [0,x^{1/p}]\times(1,\infty] \bigr).
\]
(Note that the normalization of
the $X$ variable may have to be properly scaled by $ct^p$ for some
$c>0$ to ensure the limit is a probability distribution.)

The details of this construction can be
repeated in modestly greater generality with
\eqref{eqn:defXY} modified as
%
\begin{equation}\label{eqn:defXYmod1}
(X,Y)=B(U_1,h(V_1))+(1-B)(h(U_2),V_2),
\end{equation}
with $h \in RV_p$ and $h(t)/t \to0$.
As before, $(X,Y)$ is standard regularly varying and
asymptotically independent and
%
\begin{equation}\label{eqn:bitMoreGeneral}
tP \biggl[ \biggl( \frac{X}{h(t)},\frac Yt \biggr) \in\cdot \biggr]
\stackrel{v}{\to} \mu(\cdot),
\end{equation}
where $\mu$ is given as in \eqref{eqn:defnu}. The condition $h(t)/t
\to0 $ is necessary and sufficient for $(X,Y) $ to be asymptotically
independent as can be seen by examining the calculations leading to
\eqref{eqn:asyIndep}.

\subsubsection{Mixture of independent standard regularly varying
random variables \textup{II;} negative $\rho$}\label{subsubsec:srvRV}

To exemplify case (iii) of \eqref{eqn:psi2} where $\rho<0$, suppose
\eqref{eqn:defXYmod1}, \eqref{eqn:bitMoreGeneral} still hold, $h(t)/t
\to0$ and $(U,V)$ are not asymptotically independent. Define
$\widetilde X = 1/X, \tilde h=1/h \in RV_{-p},$ and a measure
$\tilde\mu$ on $[0,\infty]\times(0,\infty]$ by
\[
\tilde\mu\bigl( [0,x]\times(y,\infty] \bigr) =\mu\biggl(
\biggl[\frac1x,\infty\biggr]\times(y,\infty] \biggr).
\]
Then
\[
tP \biggl[ \biggl( \frac{\widetilde X}{\tilde h(t)}, \frac Yt \biggr)
\in\cdot \biggr] \stackrel{v}{\to} \tilde\mu(\cdot),
\]
in $M_+ ( [0,\infty]\times(0,\infty] )$. The reason this
works is that the first space in the product $[0,\infty]\times
(0,\infty]$ is compact:
\begin{eqnarray*}
tP \biggl[ \frac{\widetilde X}{\tilde h (t)} \leq x, \frac Yt >y
\biggr]=tP \biggl[ \frac{X}{h(t)} \geq\frac1x, \frac Yt>y\biggr] \to
\mu \biggl( \biggl[\frac1x,\infty\biggr]\times(y,\infty] \biggr).
\end{eqnarray*}
So using $(\widetilde X,Y)$, we have an example of case (iii) of
\eqref{eqn:psi2} where $\rho=-p<0$, $\alpha(t)=\beta(t)=\tilde h(t).$
The conditioned limit distribution is
\[
H(x)=\lim_{t \to\infty} P[\widetilde X/\tilde h(t) \leq x|Y>t]= \mu
\biggl( \biggl[\frac1x,\infty\biggr]\times(1,\infty] \biggr).
\]

\subsubsection{Mixture of independent standard regularly varying
random variables~\textup{III;} $\rho=0$}\label{subsubsec:srvPI}
Finally, suppose \eqref{eqn:defXYmod1} still holds but this time
suppose $h \in\Pi(g) $ is nondecreasing and $\Pi$-varying with
auxiliary function $g(t)$. [E.g., we could take $h(t)=\log t, g(t)
=1.$] Then $h(t)/t \to0$ as $t \to\infty$ so $(X,Y)$ is standard
regularly varying as well as asymptotically independent. To verify this
we need the fact that if $\xi$ is either $U$ or $V$, then
%
\begin{equation}\label{eqn:cheeryfact}
tP \biggl[\frac{h(\xi)}{t}>x \biggr] \to0\qquad  (x>0, t \to\infty).
\end{equation}
To see this, let $K$ be a large number and
\begin{eqnarray*}
tP \biggl[\frac{h(\xi)}{t}>x \biggr]&=&  tP
\biggl[\frac{h(\xi)}{t}>x,\xi\leq tK \biggr]+
tP \biggl[\frac{h(\xi)}{t}>x, \xi > tK \biggr]
\\
&\leq& o(1)+tP[\xi>tK] \to K^{-1}.
\end{eqnarray*}
The upper bound is arbitrarily small and thus we verified
\eqref{eqn:cheeryfact}.

Now we check that $(X,Y)$ is standard regularly varying and
asymptotically independent:
\begin{eqnarray*}
&& tP \biggl[ \frac Xt >x  \mbox{ or } \frac Yt >y \biggr]
\\
&&\qquad =  \frac t2
P \biggl[ \frac{{U_1}}{t} >x \mbox{ or } \frac{ h({V_1})}{t}
>y \biggr]
 +\frac t2 P \biggl[ \frac{h({U_2})}{t} >x \mbox{ or }
\frac{{V_2}}{t}
>y \biggr]
\\
&&\qquad =  o(1) +\frac t2 P \biggl[ \frac{{U_1}}{t} >x \biggr] + \frac t2 P
\biggl[ \frac{ {V_2}}{t}
>y \biggr] \to \frac12 (x^{-1}+y^{-1}).
\end{eqnarray*}
Note we applied \eqref{eqn:cheeryfact}.

Next consider
\begin{eqnarray*}
&& tP \biggl[ \frac{X-h(t)}{g(t)} \leq x, \frac Yt >y\biggr]
\\
&&\qquad  =o(1) +\frac t2 P \biggl[ \frac{h(U_2)-h(t)}{g(t)} \leq x,
\frac{V_2}{t}
>y\biggr]
\\
&&\qquad \sim \frac t2 P \biggl[ \frac{U_2}{t}
\leq\frac{h^\leftarrow(g(t)x+h(t))}{t},\frac{V_2}{t}>y \biggr] \sim
\frac t2 P \biggl[ \frac{U_2}{t} \leq e^x,
\frac{V_2}{t}>y \biggr]
\\
&&\qquad \to \frac12 \nu \bigl( [0,e^x]\times(y, \infty] \bigr).
\end{eqnarray*}
This exemplifies case (ii) of \eqref{eqn:psi2} with $\rho=0, \beta
(t)=h(t) $ and $\alpha(t)=g(t)$.
The form of the conditioned limit is
\[
P \biggl[ \frac{X-h(t)}{g(t)} \leq x\Big| Y>t \biggr] \to\frac12 \nu \bigl(
[0,e^x]\times(1, \infty] \bigr)=:H(x),\qquad  x \in \mathbb{R}.
\]

\section{Change of coordinate system}\label{sec:changeCoord}
How much freedom do we have to measure the $X$-variable in different
units? This issue was raised in the discussion to \mbox{Heffernan} and Tawn
\cite{heffernantawn2004} and we try to offer further insight on the
matter here. For the example in Section \ref{subsubsec:ParetoNoNo} we
saw that for $(X,Y)$ bivariate normal, it was possible to transform
$X\mapsto\log b^\leftarrow(X) $ and get a conditional limit but the
transformation $X\mapsto b^\leftarrow(X)$ did not preserve existence of
{conditional} limits. Can something more general be said about this
issue?

Starting with \eqref{eqn:basicConvStandardForm} where the $Y$-variable
is standardized, for what monotone increasing functions $h(\cdot)$ do
there exist centering and scaling functions \mbox{$\alpha_2(t)>0$}, $\beta_2(t)
\in\mathbb{R}$, such that for some limit measure $\mu_2 $ satisfying
the nondegeneracy assumptions at the beginning of Section
\ref{sec:basicDefResults} we have
%
\begin{equation}\label{eqn:transformedStandard}
tP \biggl[ \biggl( \frac{h(X)-\beta_2(t)}{\alpha_2(t)}, \frac Yt
\biggr) \in\cdot \biggr] \stackrel{v}{\to} \mu_2
\end{equation}
in $M_+ ( [-\infty,\infty]\times(0,\infty] )$? This problem has {many}
similarities to ones considered in \cite{balkema1973,resnick1973} and
the experience gained in Section~\ref{subsubsec:ParetoNoNo} is helpful.

In \eqref{eqn:basicConvStandardForm}, assume centering by $\beta(t)$ is
really necessary; that is, suppose it is not the case that $\beta
(t)=o(\alpha(t)).$ [If $\beta (t)=o(\alpha(t))$, the following
arguments are easier and lead to regular variation of $h$.] Assume
\eqref{eqn:transformedStandard} and rewrite the left side of
\eqref{eqn:transformedStandard} evaluated on
$[-\infty,x]\times(y,\infty]$ as
\[
tP \biggl[ \frac{X-\beta(t)}{\alpha(t)} \leq
\frac{h^\leftarrow(\alpha_2(t) x +\beta_2(t)) -\beta(t)}{\alpha (t)},
\frac Yt >y \biggr].
\]
Since this converges, there must exist a limit $\psi(x)$ such
that
%
\begin{equation}\label{eqn:goToThis}
\frac{h^\leftarrow(\alpha_2(t) x +\beta_2(t)) -\beta(t)}{\alpha
(t)}
\to\psi(x)
\end{equation}
and then we see that
%
\begin{equation}\label{eqn:sharp1}
\mu\bigl( [-\infty,\psi(x)]\times(y,\infty] \bigr)= {\mu_2} \bigl(
[-\infty,x]\times(y,\infty] \bigr).
\end{equation}
The limit $\psi$ cannot be constant without violating the nondegeneracy
assumption for $\mu_2$. Inverting \eqref{eqn:goToThis} we get
\[
\frac{h ( y\alpha(t) +\beta(t) ) -\beta_2(t)}{\alpha_2(t)}
\to\psi^\leftarrow(y).
\]
This suggests we set
%
\begin{equation}\label{eqn:sharp2}
\beta_2(t)=h(\beta(t)),
\end{equation}
since
%
\begin{equation}\label{eqn:formOfh}
\frac{h ( y\alpha(t) +\beta(t) ) -h(\beta(t))}{\alpha_2(t)}
\to\psi^\leftarrow(y) -\psi^\leftarrow(0)=:\chi(y)
\end{equation}
and presuming $\chi(1)>0$, we could set
\[
\alpha_2(t)=h\bigl(\alpha(t)+\beta(t)\bigr)-h(\beta(t)).
\]

We now look at some possible forms of $h$ which allow change of
coordinate system \eqref{eqn:transformedStandard}. We do not achieve
necessary and sufficient conditions but come to an understanding of how
to generate broad classes of functions $h$ permitting nonlinear
transformation of $X$.

\subsection{Case \textup{A:} $\alpha(t) $ is asymptotically a
constant}\label{subsec:alphaConst} Assume $\beta(t) \uparrow
\infty$ as $t \to\infty$. If $\alpha\sim1$, then
\[
\frac{h ( y +\beta(t) ) -h(\beta(t))}{\alpha_2(t)}
\to\chi(y),
\]
and changing variables yields
\[
\frac{h ( y +t ) -h(t)}{\alpha_2(\beta^\leftarrow(t))}
\to\chi(y),
\]
or
%
\begin{equation}\label{eqn:sharpsharp}
\frac{h (\log tx ) -h(\log t)}{\alpha_2{ (\beta^\leftarrow(\log t) )}}
\to\chi(\log x),\qquad  x>0.
\end{equation}
Since $h\circ\log$ is nondecreasing, either \cite{gelukdehaan1987}
\begin{longlist}[(a)]
\item[(a)] $h \circ\log\in RV_p, p >0$, in which case
$\alpha_2(\beta^\leftarrow({\log t}))\sim h(\log t)$

or

\item[(b)]\label{item:b} $h \circ\log\in\Pi(\alpha_2\circ
\beta^\leftarrow(\log t)).$
\end{longlist}

Conclusion: If $\alpha\sim1$, we may change coordinates $X\mapsto
h(X)$, provided $h \circ\log\in RV_p \cup\Pi(\alpha_2\circ
\beta^\leftarrow(\log t)).$

\begin{Remark}
1. In Section \ref{subsubsec:ParetoNoNo}, $\alpha(t) =1. $ We tried
$h(x)=b^\leftarrow(x)$ but did not get a conditioned limit law. In
Section~\ref{subsubsec:ParetoNoNo}, $h\circ\log=b^\leftarrow\circ \log$
is neither regularly varying, nor $\Pi$-varying.
{\smallskipamount=0pt
\begin{enumerate}[3.]
\item[2.] In Section \ref{subsubsec:marginalsx}, $\alpha(t) =1.$ We
tried $h(x)=\log b^\leftarrow(x)$ which led to a conditioned limit law
because Lemma~\ref{lemma:Pi} proved $h\circ\log= \log
b^\leftarrow\circ\log\in\Pi(\log).$

\item[3.] The result in (b) suggests how to construct other examples of
$h$ which lead to conditioned limits. If $g$ is any slowly varying
function, then $\int_1^x g(u)u^{-1}du $ is $\Pi$-varying with auxiliary
function $g$ (\cite{dehaan1976}, \cite{resnick1987}, page~30). Define $h$ by
$h(\log x)=\int_1^x g(u)/u \,du $ or
\[
h'(x)=g(e^x),\qquad  h(x)=\int_0^x g(e^u)\,du.
\]
Any such $h$ will lead to a conditioned limit.
Examples include:
\begin{itemize}
\item $g(x) =\log x$ and $h(x) =x^2/2.$

\item$g(x)=\log\log x$ and $h(x)=\int_0^x \log u\, du \sim x\log x.$

\item$g(x)=(\log x)^p $ and $h(x)=\frac{x^{p+1}}{p+1}$ for $p>0$.
\end{itemize}
For an example where $h\circ\log\in
RV_p$ for $p>0$, set
\[
h(\log x)=U(x)\in RV_p\quad\mbox{or}\quad h(x)=U(e^{ x}).
\]
Apply this to the convergence \eqref{eqn:xyStandard} for the bivariate
normal pair $(X,Y)$ where recall
\[
\beta(t) =\rho b(t),\qquad \alpha(t)=1,\qquad \mu\bigl(
[-\infty,x]\times(y,\infty] \bigr) =N \biggl(\frac{x}{\sqrt{1-\rho^2}}
\biggr) y^{-1}.
\]
Then evaluating \eqref{eqn:sharpsharp} with $h(\log t)=U(t) \in RV_p,
p>0,$ gives, with $\alpha_2\circ\beta^\leftarrow {\circ\log}=U$ that
\[
\frac{U(tx)-U(t)}{U(t)} \to x^p -1=\chi(\log x).
\]
Therefore, $\chi(y)=e^{py} -1,$ and from \eqref{eqn:sharp1}
\begin{eqnarray*}
tP \biggl[ \frac{U(e^X) -U(e^{\rho b(t)})}{U(e^{\rho b(t)})} \leq
x,\frac{b^\leftarrow(Y)}{t} >y \biggr]&=&
\mu\bigl([-\infty,\chi^\leftarrow(x)]\times(y,\infty] \bigr)
\\
&=& N \biggl( \frac{p^{-1}\log(1+x)}{\sqrt{1-\rho^2}} \biggr) y^{-1}.
\end{eqnarray*}

So for this example, $\beta_2 (t)=\alpha_2(t)=U(e^{\rho b(t)}).$
\end{enumerate}}
\end{Remark}

\subsection{Case \textup{B:} $\alpha(t)$ is not asymptotically a constant}\label{subsec:alphaNotConst}
Again assume $\beta(t) \uparrow \infty$ as $t \to\infty$.
Transform~\eqref{eqn:formOfh} to get
%
\begin{equation}\label{eqn:formOfh2}
\frac{h ( y\alpha\circ\beta^\leftarrow(t) +t )
-h(t)}{\alpha_2\circ\beta^\leftarrow(t)}
\to\chi(y)
\end{equation}
which is of the form
\[
\frac{h ( t+f(t)y )
-h(t)}{\alpha_*(t)} \to\chi(y).
\]
To proceed further in a way that generates a broad class examples,
suppose $f(t)=\alpha\circ\beta^\leftarrow(t)$ is
\textit{self-neglecting} \cite{binghamgoldieteugels1987}. A simple
sufficient condition is $f'(t) \to0$ and $f$ self-neglecting means it
is the auxiliary function of a {$\Gamma$}-varying function (see
Appendix~\ref{subsubsec:fctClass}) and that
\[
H(x):=\exp\biggl\{\int_1^x \frac{1}{f(u)}\,du \biggr\} \in\Gamma(f).
\]
Then defining the function $V$ by
\[
h=V\circ H \mbox{ or equivalently } V=h\circ H^\leftarrow
\]
we have either (\cite{dehaan1976}, page~249, \cite{resnick1987},
page~36)
\begin{longlist}[(a)]
\item[(a)] $V \in\Pi$ and $\chi(y)=\log e^y =y$;
\end{longlist}
or
\begin{longlist}[(a)]
\item[(b)] $V \in RV_p, p>0$ and $\chi(y)=e^{py}-1.$
\end{longlist}

Conclusion: We considered the case that $\beta\neq o(\alpha)$ and
$\beta(t) \uparrow\infty$ and $\alpha$ not asymptotically a constant.
For such a case, the change of variable $X\mapsto h(X)$ preserves
conditioned limits provided $h$ is either the composition of a
\mbox{$\Pi$-}varying function and a $\Gamma$-varying function or the
composition of a regularly varying function and a \mbox{$\Gamma$-}varying
function. (The composition of a regularly varying function and a
\mbox{$\Gamma$-}varying function is another $\Gamma$-varying function; see
\cite{dehaan1970}, \cite{resnick1987}, page~36).

\section{Discussion and concluding remarks}\label{sect:discuss}
The statistical models proposed by \mbox{Heffernan} and Tawn
\cite{heffernantawn2004} are based on the assumption that for $(X,Y)$
having Gumbel marginal distributions, there exist normalizing functions
$\alpha(\cdot)$ and $\beta(\cdot)$ such that the conditional
distribution of $(X-\beta(y))/\alpha(y)$ given $Y=y$ can be
approximated for large $y$ by some nondegenerate, proper $G(x)$. We
have built our theory by standardizing $Y$ to have asymptotically
Pareto distribution and looked at the conditional distribution of
$(X-\beta(t))/\alpha(t)$ given $Y>t$ which also leads to conditional
distributions for $(X-\beta(Y))/\alpha(Y)$ given $Y>t$. This
formulation is consistent with the \mbox{Heffernan} and Tawn
\cite{heffernantawn2004} approach and allows a mathematically precise
theory which can be related to the extended theory of multivariate
regular variation.

From the perspective of statistical modeling, important results
are contained in Propositions~\ref{prop:randomNorming}
and~\ref{prop:whoKnows}. These propositions reveal the
factorization of the limit distribution obtained when $X$ is
normalized by the value of $Y$ that occurs with it. This
factorization permits a significant simplification of models based
on the limit form, as it enables the assumption of limiting
independence between the conditioning and standardized variables.
This independence assumption was employed in~\cite{heffernantawn2004}
and is key to statistical modeling and extrapolation.

One issue we have not resolved is consistency of different models.
The definition \eqref{eqn:basicConv} or its standardized version
\eqref{eqn:basicConvStandardForm} is not symmetric in the $X,Y$
variables. However, when fitting models to data one has a choice of
which variable to condition being large and a logical issue is
whether the various models obtained by conditioning on different
variables are related to each other in any way. Conditions for
consistency would strengthen the statistical model assumptions based
on this representation and therefore potentially improve the
ability of such approaches to describe the joint distribution in tail
regions where there is naturally little data. Currently we have
nothing terribly useful to say on this issue other than to point out
that it seems important to understand consistency better.

\begin{appendix}
\section*{Appendices}\label{sec:appendices}
For convenience, this section collects some notation, needed
background on regular variation
and notions on vague convergence needed for some formulations and proofs.

\renewcommand{\thesubsection}{A.\arabic{subsection}}
\setcounter{subsection}{0}

\subsection{Vector notation}\label{subsubsec:vecNotation}
Vectors are denoted by bold letters, capitals for random vectors and
lower case for nonrandom vectors. For example:
$\bx=(x^{(1)},\ldots,x^{(d)}) \in\mathbb{R}^d.$ Operations between
vectors should be interpreted componentwize so that for two vectors
$\bx$ and $\bz$
\begin{eqnarray*}
\bx &<& \bz \mbox{ means }x^{(i)}<z^{(i)},\qquad i=1,\dots,d,
\\
\bx &\leq& \bz \mbox{ means }x^{(i)}\leq z^{(i)},\qquad i=1,\dots,d,
\\
\bx&=&  \bz \mbox{ means }x^{(i)}=z^{(i)},\qquad i=1,\dots,d,
\\
\bz\bx&=& \bigl(z^{(1)} x^{(1)},\dots,z^{(d)} x^{(d)}\bigr),
\\
\bx\vee\bz&=&  \bigl(x^{(1)} \vee z^{(1)},\dots,x^{(d)} \vee z^{(d)}\bigr),\qquad
\frac{\bx}{\bz} = \biggl( \frac{x^{(1)} }{z^{(1)} },\dots,\frac{x^{(d)}
}{z^{(d)} } \biggr),
\end{eqnarray*}
and so on. Also define $\mathbf 0= (0,\dots,0).$ For a real
number $c$, denote as usual $c\bx=(cx^{(1)},\dots,cx^{(d)}).$ We
denote the rectangles (or the higher dimensional intervals) by
\[
[\mathbf a,\mathbf b]=\{\bx\in{\mathbb{R}^d}\dvtx  \mathbf a
\leq\bx\leq\mathbf b\}.
\]
Higher dimensional rectangles
with one or both endpoints open are defined analogously, for
example,
\[
(\mathbf a,\mathbf b]=\{\bx\in{\mathbb{R}^d}\dvtx
\mathbf a < \bx\leq\mathbf b\}.
\]

\subsection{The function classes $\Pi$ and $\Gamma$}\label{subsubsec:fctClass}
Continue the domain of attraction discussion:
Writing \eqref{eqn:doa'} as
\[
\biggl(\frac{1}{1-F(a(t)x+b(t))} \biggr)\Big/t \to(1+\gamma x )^{1/\gamma}
\]
and inverting yields as $t \to\infty$
%
\begin{equation}\label{eqn:PiEtc}
\frac{b(ty)-b(t)}{a(t)} \to
\cases{
\dfrac{y^\gamma-1}{\gamma},&\quad if $\gamma\neq0$,
\cr\noalign{}
\log y,&\quad if $\gamma=0$.}
\end{equation}
In case $\gamma=0$, \eqref{eqn:PiEtc} says that $b(\cdot) \in\Pi
(a (\cdot) )$; that is, the function $b(\cdot)$ is
$\Pi$-varying with auxiliary function $a(\cdot)$
(\cite{resnick1987}, pages 26ff, \cite{binghamgoldieteugels1987,gelukdehaan1987,dehaan1970}).

More generally (\cite{binghamgoldieteugels1987},
Chapter~3, \cite{dehaanresnick1979}) define for an auxiliary function $a(t)>0$,
$\Pi_+(a)$ to be the set of all functions $\pi\dvtx \mathbb{R}_+ \mapsto
\mathbb{R}_+$ such that
%
\begin{equation}\label{eqn:piplus}
\lim_{t \to\infty}\frac{\pi(tx)-\pi(t)}{a(t)}=k\log x,\qquad  x>0,
k>0.
\end{equation}
The class $\Pi_- (a)$ is defined similarly except that $k<0$ and
\[
\Pi(a)=\Pi_+ (a)\cup\Pi_- (a).
\]
By adjusting the auxiliary function in the denominator, it is always
possible to assume $k=\pm1$.

Two functions $\pi_i \in\Pi_\pm(a) $, $i=1,2$, are $\Pi(a)$-equivalent if for some
$c\in\mathbb{R}$
\[
\lim_{t \to\infty}\frac{\pi_1(t)-\pi_2(t)}{a(t)}= c.
\]
There is usually no loss of generality in assuming $c=0$.

The class of regularly varying functions with index $\rho\in
\mathbb{R}$ is denoted by $RV_\rho$ so that $U\dvtx \mathbb{R}_+\mapsto
\mathbb{R}_+$ satisfies $U \in RV_\rho$ if
%
\begin{equation}\label{eqn:defRV}
\lim_{t \to\infty}\frac{U(tx)}{U(t)}=x^\rho,\qquad  x>0.
\end{equation}

The following are known facts about $\Pi$-varying functions.
\begin{enumerate}
\item\label{item:fact1} We have $\pi \in\Pi_+(a)$ iff $1/\pi\in\Pi
_-(a/\pi^2).$

\item\label{item:fact2} If $\pi \in\Pi_+(a)$, then (\cite{binghamgoldieteugels1987}, page~159 or
\cite{dehaanresnick1979}, page~1031)
there exists a continuous and strictly increasing $\Pi(a)$-equivalent
function $\pi_0$ with $\pi-\pi_0= o(a)$.

\item\label{item:fact3} If $\pi \in\Pi_+(a)$, then
\[
\lim_{t \to\infty} \pi(t)=:\pi(\infty)
\]
exists. If $\pi(\infty)=\infty$, then $\pi\in RV_0$ and $\pi(t)/a(t)
\to\infty$. If $\pi(\infty)<\infty$, then $\pi(\infty)-\pi(t) \in
\Pi_- (a)$ and $\pi(\infty)-\pi(t) \in RV_0$ and
$(\pi(\infty)-\pi(t))/\break a(t) \to\infty$. (Cf. \cite{gelukdehaan1987}, page~25.) Furthermore,
\[
\frac{1}{\pi(\infty)-\pi(t)} \in\Pi_+ \bigl(a/\bigl(\pi(\infty)-\pi(t)\bigr)^2 \bigr).
\]
\end{enumerate}

In addition to the function class $\Pi$ we need de Haan's class
$\Gamma$ (\cite{binghamgoldieteugels1987,gelukdehaan1987,dehaan1970,dehaan1974,resnick1987}). A function
$V\dvtx \mathbb{R}_+ \mapsto\mathbb{R}_+$ is a $\Gamma$-function with
auxiliary function $f$ [written $V \in\Gamma(f)$] if, as $t \to
\infty,$
\[
\frac{V(t+xf(t))}{V(t)} \to e^x,\qquad  x>0.
\]
For $V$ nondecreasing, $V \in\Gamma(f)$ iff $V^\leftarrow\in\Pi(f\circ
V^\leftarrow).$

\subsection{Vague convergence}\label{subsubsec:vagueConv}
For a nice space $\mathbb{E}$, that is, a space which is locally compact with
countable base (e.g., a finite dimensional Euclidean space),
denote $M_+(\mathbb{E})$ for the nonnegative Radon measures on Borel subsets of
$\mathbb{E}$. This space is metrized by the vague metric. The notion of vague
convergence in this space is as follows: If $\mu_n \in M_+(\mathbb{E})$ for $n
\geq0$, then $\mu_n$ converge vaguely to $\mu_0$ (written $\mu_n
\stackrel{v}{\to} \mu_0$) if for all bounded continuous functions $f$
with compact support we have
\[
\int_{\mathbb{E}} f\,d\mu_n \to
\int_{\mathbb{E}} f\,d\mu_0\qquad (n \to\infty).
\]
This concept allows us to write \eqref{eqn:doa'} as
%
\begin{equation}\label{eqn:doa''}
tP\biggl[\frac{Y-b(t)}{a(t)} \in\cdot\biggr] \stackrel{v}{\to} m_\gamma(\cdot),
\end{equation}
vaguely in $M_+((-\infty,\infty])$ where
\[
m_\gamma((x,\infty])=(1+\gamma x)^{-1/\gamma}.
\]
Standard references include \cite{kallenberg1983,neveu1977} and \cite{resnick1987},
Chapter~3.
\end{appendix}

\section*{Acknowledgments}
The authors would like to thank two anonymous referees for careful
reading of the initial submission, and their comments which contributed
to a substantially improved final version.

Thanks to Lancaster University and Cornell University's School of
Operations Research and Industrial Engineering and Department of
Statistics for funding and hospitality during a visit to Cornell in
September 2004.

\printaddresses

\end{document}